\documentclass[10pt]{article}

\usepackage{amssymb,latexsym}
\usepackage{epsfig}
\usepackage{eufrak}
\usepackage{amsmath}
\usepackage{mathrsfs}
\usepackage{color}

\newtheorem{theorem}{Theorem}
\newtheorem{lemma}{Lemma}
\newtheorem{corollary}{Corollary}

\newcommand{\be}{\begin{equation}}
\newcommand{\ee}{\end{equation}}
\newcommand{\bee}{\begin{eqnarray*}}
\newcommand{\eee}{\end{eqnarray*}}
\newcommand{\bel}{\begin{eqnarray}}
\newcommand{\eel}{\end{eqnarray}}
\newcommand{\bec}{\begin{cases}}
\newcommand{\eec}{\end{cases}}
\newcommand{\bem}{\begin{bmatrix}}
\newcommand{\eem}{\end{bmatrix}}

\newcommand{\la}{\label}
\newcommand{\li}{\left}
\newcommand{\ri}{\right}

\newcommand{\ovl}{\overline}

\newcommand{\ep}{\epsilon}
\newcommand{\vep}{\varepsilon}
\newcommand{\lm}{\lambda}

\newcommand{\si}{\sigma}

\newcommand{\de}{\delta}

\newcommand{\vDe}{\varDelta}

\newcommand{\se}{\theta}

\newcommand{\ze}{\zeta}
\newcommand{\al}{\alpha}
\newcommand{\ba}{\beta}

\newcommand{\ro}{\rho}
\newcommand{\ka}{\kappa}
\newcommand{\om}{\omega}
\newcommand{\Om}{\Omega}

\newcommand{\f}{\frac}
\newcommand{\sq}{\sqrt}
\newcommand{\cd}{\cdots}

\newcommand{\qu}{\quad}
\newcommand{\qqu}{\qquad}

\newcommand{\mscr}{\mathscr}
\newcommand{\mcal}{\mathcal}

\newcommand{\bb}{\mathbb}
\newcommand{\fra}{\mathfrak}

\newcommand{\bs}{\boldsymbol}

\newcommand{\tx}{\text}

\newcommand{\iy}{\infty}

\newcommand{\pa}{\partial}

\newcommand{\bed}{\begin{description}}
\newcommand{\eed}{\end{description}}
\newcommand{\bei}{\begin{itemize}}
\newcommand{\eei}{\end{itemize}}
\newcommand{\ben}{\begin{enumerate}}
\newcommand{\een}{\end{enumerate}}
\newcommand{\bib}{\bibitem}
\newcommand{\beL}{\begin{lemma}}
\newcommand{\eeL}{\end{lemma}}
\newcommand{\beT}{\begin{theorem}}
\newcommand{\eeT}{\end{theorem}}
\newcommand{\beC}{\begin{corollary}}
\newcommand{\eeC}{\end{corollary}}
\newcommand{\sect}{\section}

\newcommand{\bpf}{\begin{pf}}
\newcommand{\epf}{\end{pf}}
\newcommand{\bsk}{\bigskip}

\newcommand{\pfbox}{\hfill\mbox{$\Box$}}

\newenvironment{pf}{\paragraph*{Proof{\rm.}}}{\pfbox\bigskip}

\begin{document}

\title{{\bf A Geometric Approach for Bounding Average Stopping Time} }

\author{Xinjia Chen\\
Department of Engineering  Technology\\
 Northwestern State
University, Natchitoches, LA 71497\\
Email: chenx@nsula.edu \quad Tel: (318)357-5521  \quad Fax:
(318)357-6145}

\date{ }

\maketitle

\begin{abstract}

We propose a geometric approach for bounding average stopping times
for stopped random walks in discrete and continuous time. We
consider stopping times in the hyperspace of time indexes and
stochastic processes. Our techniques relies on exploring geometric
properties of continuity or stopping regions. Especially, we make
use of the concepts of convex sets and supporting hyperplane.
Explicit formulae and efficiently computable bounds are obtained for
average stopping times.  Our techniques can be applied to bound
average stopping times involving random vectors, nonlinear stopping
boundary, and constraints of time indexes. Moreover, we establish a
stochastic characteristic of convex sets and generalize Jensen's
inequality, Wald's equations and Lorden's inequality, which are
useful for investigating average stopping times.

\end{abstract}

%\tableofcontents

\sect{Introduction}

\bsk

In many areas of engineering and sciences, especially probability
and statistics, it is interested to investigate the expectation of
stopping times defined in terms of stochastic processes of
stationary and independent increments. For example, a frequent topic
of random walk \cite{Berg, Gut2009, Lawler, Rudnick, Spitzer}
concerns a stopping time which is the smallest positive integer $n$
such that the partial sum $X_1 + \cd + X_n$ is greater than $f(n)$,
where $X_1, X_2, \cd$ are i.i.d. random variables and $f$ is a
function of $n$. Since many sequential hypothesis testing and
estimation procedures can be cast into the context of such stopping
time, for analyzing the efficiency of statistical inference, it is
 of practical importance to evaluate the expectation of  such stopping time in the area of sequential analysis
 \cite{BKGosh, SEqGhosh, Govindarajulu, Lai, MukSi, Tartakovsky, Wald}.  Although the
 literature on such stopping time is abundant, most existing works are focused on the asymptotic analysis of average stopping times
 (see, e.g., \cite{Siegmund, Woodroofe} and the references therein).
 Existing techniques such as Lorden's inequality \cite{Lorden}
 for bounding average stopping times are limited to very specific forms of $f(n)$.  In many practical situations,
 $f(n)$ can be complicated functions without nice properties such as linearity and monotonicity. The time index $n$ may be
 restricted to a subset of natural numbers, as usually required in group sequential methods
 \cite{Bartroff, Jennsion, Mukhopadhyay, Proschan, Whitehead}. The underlying variables $X_i$ may be random
 vectors.  However, there lacks of effective technique for
 obtaining tight bounds for average stopping times, which are general enough to
 deal with the nonlinearity of the function $f(n)$, the constraint of the time index $n$,
  and the dimensionality of random variables $X_1, X_2,
 \cd$.  Motivated by this situation,  we propose a geometric approach to bound  average stopping times in a general setting. We consider
 stopping times in the hyperspace of the tuple $(n, X_1+ \cd + X_n)$, where $X_i$ are allowed to be random vectors and $n$
 is contained by a subset $\mscr{N}$ of natural numbers.
 A stopping time is represented as the first time $n \in \mscr{N}$ that the tuple $(n, X_1+ \cd + X_n)$
 falls into a certain region, referred to as a stopping region (or equivalently, falls outside of a certain region, referred to as a continuity region).
Our main idea is to make use of the geometric properties of the continuity region or stopping region. Particularly, we will use concepts such as
convexity and supporting hyperplane to develop bounds for average stopping times, which are either explicit or amenable for convex minimization.

The remainder of the paper is organized as follows. In Section 2, we
propose to investigate stopping times in a geometric setting, which
makes it possible to use geometric concepts such as convex hull,
convex set, and supporting hyperplane, etc.   Afterward, we
establish a probabilistic property of convex sets, which plays a
crucial role in bounding average stopping times.  In Section 3, we
generalize Jensen's inequality, Wald's equations and Lorden's
inequality, which are fundamental tools for investigating average
stopping times.  In Section 4, we establish efficient convex
minimization techniques for bounding average stopping times.  In
Section 5, we develop explicit formulae for bounding average
stopping times by virtue of the concept of supporting hyperplane. In
Section 6, we propose to bound average stopping times by combining
the power of concentration inequalities and the concept of geometric
convexity. In Section 7, we extend the techniques to bound average
stopping times relevant to L\'{e}vy processes.  Section 8 is the
conclusion. Most proofs are given in Appendices.

In this paper, we shall use the following notations.  An empty set
is denoted by $\emptyset$. The infimum of an empty set is defined as
$\iy$. The supremum of an empty set is defined as $0$. For a set
$\mscr{S}$, its closure and boundary are denoted by $\ovl{\mscr{S}}$
and $\pa \mscr{S}$, respectively.

The set of positive integers is denoted by $\bb{N}$. The set of
non-negative integers is denoted by $\bb{Z}^+$. The set of real
numbers is denoted by $\bb{R}$. The set of non-negative real numbers
is denoted by $\bb{R}^+$. The set of real-valued column matrices of
size $d \times 1$ is denoted by $\bb{R}^d$.  A column matrix in
$\bb{R}^d$ is also called a vector. The notation $\bs{0}_d$ denotes
a column matrix of size $d \times 1$ with all elements being $0$.
The notation $\bs{1}_d$ denotes a column matrix of size $d \times 1$
with all elements being $1$.

We use notation $\top$ to denote the transpose of a matrix.  We
define the following operations of column (or row) matrices:

$A B$ denotes  the product of {\small $A = [a_1, \cd, a_d]^\top$}
and {\small $B = [b_1, \cd, b_d]^\top$} in the sense that {\small $A
B = [ a_1 b_1, \cd, a_d b_d ]^\top$}.

$\f{A} {B}$ denotes the quotient  of $A =[a_1, \cd, a_d]^\top$
divided by $B = [b_1, \cd, b_d]^\top$ in the sense that
 {\small $\f{A}{B} = [ \f{a_1}{b_1}, \cd, \f{a_d}{b_d} ]^\top$}.

For $A = [a_1, \cd, a_d]^\top$, we use $A^i$ to denote the $i$-th
power of $A$ in the sense that $A^i = [a_1^i, \cd, a_d^i]^\top$.
Similarly, for $A = [a_1, \cd, a_d]$, we use $A^i$ to denote the
$i$-th power of $A$ in the sense that $A^i = [a_1^i, \cd, a_d^i]$.

For $A = [a_1, \cd, a_d]^\top$, we use $|A|$ to denote the absolute
value of $A$ in the sense that $|A| = [|a_1|, \cd, |a_d|]^\top$.
Similarly, for $A = [a_1, \cd, a_d]$, we use $|A|$ to denote the
absolute value of $A$ in the sense that $|A| = [|a_1|, \cd, |a_d|]$.

For matrices $A = [a_1, \cd, a_d]^\top$ and
 $B = [b_1, \cd, b_d]^\top$, we write $A \leq B$ if $a_i \leq b_i$ for $i = 1, \cd, d$.

The Euclidean norm of a column matrix or row matrix is denoted by
$||.||$.

For a function,  $f(v)$, of $v = [v_1, \cd, v_d]^\top \in \bb{R}^d$,
we use $\f{\pa f(v) } {\pa v}$ to denote the gradient of $f(v)$ with
respect to $v$, that is, {\small $\f{\pa f(v) } {\pa v} = \li [
\f{\pa f(v) } {\pa v_1}, \cd, \f{\pa f(v) } {\pa v_d} \ri ]$}.

The probability space is denoted by $(\Om, \mscr{F}, \Pr)$, where
$\Om$ is the sample space, $\mscr{F}$ is the $\si$-algebra on $\Om$,
and $\Pr$ is the probability measure.  The probability of an event
$E$ is denoted by $\Pr \{E\}$.   The mathematical expectation of a
random variable (scalar or vector) $X$ is denoted by $\bb{E} [ X]$.
 Let $\bb{I}_E$ denote the indicator function such that it assumes value $1$ if the event $E$ occurs and
 it assumes value $0$ otherwise.

For random vector $X = [ \bs{x}_1, \cd, \bs{x}_d ]^\top$, we define
$X^+ = [ \max(0, \bs{x}_1), \cd, \max(0, \bs{x}_d)]^\top$ as the
non-negative part of $X$.  Similarly, we define $X^- = [ \max(0, -
\bs{x}_1), \cd, \max(0, - \bs{x}_d)]^\top$ as the non-positive part
of $X$.

 The other notations will be made clear as we proceed.

\section{Stopping Times and Convex Sets}

In this section, we shall propose to investigate stopping times with their geometric representations.  We shall also establish a connection
between stopping times and convex sets.  A stochastic characterization of convex sets is developed.

\subsection{Geometric Representation of Stopping Time}
\la{sec21sub}

Existing methods for bounding the average of a stopping time
typically focus on exploring the properties of the function defining
the stopping time. Consider, for example, the stopping time
mentioned in the introduction of this paper.  To bound the
expectation of stopping time \be \la{st8899} \bs{N} = \inf \{n \in
\bb{N}: X_1 + \cd + X_n > f(n) \}, \ee
 conventional wisdom is to explore the function $f(n)$ for properties such as
linearity and  monotonicity  which could be useful for bounding the
average stopping time.  We would like to point out that the methods
in this direction usually fail to fully exploit the geometric
information of the underlying continuity or stopping regions.  To
clearly address this point, we shall first provide geometric
representations of stopping times in the sequel.

Throughout the remainder of this paper, we shall use the following
notations and definitions.  Let $0 \leq N_0 < N_1 < N_2 < \cd $ be
an increasing sequence of integers and define $\mscr{N} = \{ N_1,
N_2, \cd \}$. Let $\mscr{R}$ be a closed subset of $\{ (t, s): t \in
\bb{R}^+, \; s \in \bb{R}^d \}$ which contains $(0, \bs{0}_d)$. Let
$X = [\bs{x}_1, \cd, \bs{x}_d ]^\top$ be a $d$-dimensional
real-valued random vector with mean $\mu = \bb{E} [ X ]$. Let $X_1,
X_2, \cd$ be i.i.d. random vectors having the same distribution as
$X$. Define $S_0 = 0$ and
\[
S_n = \sum_{i=1}^n X_i, \qqu \ovl{X}_n = \f{S_n}{n}
\]
for $n \in \bb{N}$.  Our effort will be devoted to stopping times
which are defined in terms of the partial sum $S_n$ (or
equivalently, empirical mean $\ovl{X}_n$), the region $\mscr{R}$ and
the set $\mscr{N}$. The stopping times defined in this way can be
fairy general.

A stopping time can be defined in terms of $S_n$ as \be \la{defFPT}
\bs{N} = \inf \{ n \in \mscr{N}: (n, S_n) \notin \mscr{R} \}. \ee We
call this expression a geometric representation of stopping time,
since it is regarding the inclusion of a random point  $(n, S_n)$ by
a domain in the Euclidean space.    The stopping time in
(\ref{defFPT}) is associated with the stopping rule: Continue
observing $S_n$ until $(n, S_n) \notin \mscr{R}$ for some $n \in
\mscr{N}$. In probabilistic terminology, $\{ S_n \}$ is called a
{\it random walk}, and the stopping time $\bs{N}$ is also called the
{\it first passage time} (FPT). Clearly, the support of the FPT is
$\mscr{N}$. For such stopping rule, the region $\mscr{R}$ is
referred to as a continuity region. The complement of $\mscr{R}$,
denoted by $\mscr{R}^c$,  is called a stopping region.

Despite the generality of the above geometric representation,
stopping times are conventionally expressed in algebraic forms. A
familiar example is the stopping time defined by (\ref{st8899}).  In
this paper, we propose to investigate stopping times based on their
geometric representations. The primary reason is that the bounding
of average stopping times can be much more easier by exploiting the
geometric properties of the underlying continuity or stopping
region. As will be seen later, this is especially true when the
continuity region or stopping region is convex. We discovered that,
for a wide variety of stopping times in the context of sequential
hypothesis testing and estimation,  the corresponding continuity or
stopping regions in geometric representations are actually convex.
In the worse case that the continuity or stopping regions  are not
convex, it is still possible to bound the average stopping time by
replacing them with their convex hulls, at the price of extra
conservatism.

To illustrate the advantage of geometric representations, consider
stopping time \[ \bs{N} = \inf \{ n \in \bb{N}: f(n, S_n) > 0 \},
\]
where $f(t, s)$ is a bivariate function of $t \in \bb{R}^+$ and $s
\in \bb{R}$.  Clearly, the continuity region is \be \la{contreg}
\mscr{R} = \{ (t, s): t \in \bb{R}^+, \; s \in \bb{R}, \;  f(t, s)
\leq 0 \} \ee and the stopping time $\bs{N} = \inf \{ n \in \bb{N}:
(n, S_n) \notin \mscr{R} \}$.  Similarly, the stopping region is \be
\la{stopreg} \mscr{R}^c = \{ (t, s): t \in \bb{R}^+, \; s \in
\bb{R}, \; f (t, s) > 0\} \ee and  the stopping time $\bs{N} = \inf
\{ n \in \bb{N}: (n, S_n) \in \mscr{R}^c \}$.

It can be shown that if $f$ is a convex function, then the
continuity region (\ref{contreg}) is convex. If $f$ is a concave
function, then the stopping region (\ref{stopreg}) is convex.  It is
important to note that the convexity of the stopping or continuity
region may also hold in situations when the function $f$ is neither
convex nor concave. Moreover, even if neither the continuity region
nor the stopping region is convex, we may still be able to bound the
average stopping time by using their convex hulls. This example
demonstrates that, in contrast to using algebraic forms of stopping
times, it is possible to exploit the convexity of the continuity or
stopping regions in geometric representations under much weaker
conditions.

\subsection{ A Stochastic Characteristic of Convex Sets}

As discussed in Section \ref{sec21sub}, there exists a useful
connection between stopping times and convex sets.   Since
continuity or stopping regions are convex in many situations, it is
natural to consider the question of under what conditions the
expectation of a random vector will be contained by a convex set.
Our investigation indicates that if a set in a finite-dimensional
Euclidean space is convex, then the set contains the expectation of
any random vector almost surely contained by the set.  More
formally, we have established the following result.

\beT

\la{fundamental}

If $\mscr{D}$ is a convex set in $\bb{R}^n$, then $\bb{E} [ \bs{\mcal{X}} ] \in \mscr{D}$ holds
for any random vector $\bs{\mcal{X}}$ such that $\Pr \{ \bs{\mcal{X}}
\in \mscr{D} \} = 1$ and that $\bb{E} [ \bs{\mcal{X}} ] $ exists.
  \eeT

See Appendix \ref{fundamentalapp} for a proof.  The converse of
Theorem  \ref{fundamental} asserts that if $\mscr{D}$ is a set in
$\bb{R}^n$ such that $\bb{E} [ \bs{\mcal{X}} ] \in \mscr{D}$ holds
for any random vector $\bs{\mcal{X}}$ such that $\Pr \{
\bs{\mcal{X}}  \in \mscr{D} \} = 1$ and that $\bb{E} [ \bs{\mcal{X}}
] $ exists, then $\mscr{D}$ is convex. This assertions is well known
and is a direct consequence of the definition of a convex set.

Theorem \ref{fundamental} plays a fundamental role in our approach
for bounding average stopping times. Moreover,  Theorem
\ref{fundamental} immediately implies Jensen's inequality. To see
this, note that if a function is convex, then its epigraph, the
region above its graph, is a convex set. Hence, if $f$ is a convex
function, then for any random variable $X$, since $(X, \; f(X))$ is
contained by the epigraph of $f$, it follows from Theorem
\ref{fundamental} that $(\bb{E}[ X ], \; \bb{E} [ f(X) ] )$ is
contained by its epigraph. This implies that $\bb{E} [ f(X) ] \geq f
(\bb{E} [ X ] )$ by the notion of epigraph.

\section{Generalizations of Jensen's Inequality, Wald's Equations and Lorden's Inequality}

In this section, we shall generalize Jensen's inequality, Wald's equations and Lorden's inequality, which can be useful for evaluating average
stopping times.

\subsection{Generalization of Jensen's Inequality}

We have derived the following results.

\beT

\la{genJen} Let $\bs{Z}$ be a random vector and $Y$ be a scalar
random variable such that $\f{ \bs{Z}  }{Y}$  and $\f{ \bb{E} [
\bs{Z} ] }{ \bb{E} [ Y ] } $ are contained by a convex set
$\mscr{D}$ in $\bb{R}^n$. Assume that $g(z)$ is a convex function of
$z \in \mscr{D}$. Then,
\[
\bb{E} \li [ Y g \li ( \f{\bs{Z}}{Y} \ri ) \ri ] \geq  \bb{E} [ Y ] g \li ( \f{ \bb{E} [ \bs{Z} ]}{\bb{E} [ Y ]} \ri ) \qqu \tx{if $Y$ is a
positive random variable};
\]
\[
\bb{E} \li [ Y g \li ( \f{\bs{Z}}{Y} \ri ) \ri ] \leq  \bb{E} [ Y ] g \li ( \f{ \bb{E} [ \bs{Z} ]}{\bb{E} [ Y ]} \ri ) \qqu \tx{if $Y$ is a
negative random variable}.
\]

\eeT

\bsk

See Appendix \ref{genJenapp} for a proof.   It should be noted that Theorem \ref{genJen} generalizes Jensen's inequality.  In the special case
that $Y = 1$, the first inequality of Theorem \ref{genJen} reduces to Jensen's inequality.

\subsection{Generalization of Wald's Equations}

Making use of Theorem \ref{genJen}, we have generalized Wald's equations \cite{Wald} as follows.

\beT

\la{genWald}

Let $X_1, X_2, \cd$ be i.i.d. random vectors having the same
distribution as $X$ with mean $\mu = \bb{E} [ X ]$ and variance $\nu
= \bb{E} [ | X - \mu |^2 ]$. Assume that $\bs{N}$ is a positive
integer-valued  random variable such that $\bb{E} [ \bs{N} ] < \iy$
and that for any possible value $n$ of $\bs{N}$, the event $\{
\bs{N} = n \}$ depends only on $X_1, \cd, X_n$.  Define $S_{\bs{N}}
= \sum_{i=1}^{\bs{N}} X_i, \; \ovl{X}_{\bs{N}}  = \f{S_{\bs{N}}
}{\bs{N}}$ and $\ovl{V}_{ \bs{N} } = \f{(S_{\bs{N}} - \bs{N} \mu)^2
}{\bs{N}}$.  The following assertions hold.

(I): If $g$ is a convex function on a convex set $\mscr{D}$ in
$\bb{R}^d$ such that $\mscr{D}$ contains $\mu$ and the range of
$\ovl{X}_{\bs{N}} $, then $\bb{ E } \li [ \bs{N} g (
\ovl{X}_{\bs{N}} ) \ri ] \geq  \bb{ E } [ \bs{N} ] g (\mu)$.

(II): If $g$ is a convex function on a convex set $\mscr{D}$ in
$\bb{R}^d$ such that $\mscr{D}$ contains $\nu$ and the range of
$\ovl{V}_{ \bs{N} }$, then $\bb{ E } \li [ \bs{N} g \li ( \ovl{V}_{
\bs{N} } \ri  ) \ri ] \geq  \bb{ E } [ \bs{N} ] g (\nu)$.

\eeT

See Appendix \ref{genWaldapp} for a proof.   To see why the inequality in the first assertion of Theorem \ref{genWald} is a generalization of
Wald's first equation, consider function $g(x) = x$. By the convexity of $g(x)$, we have $\bb{ E } \li [ \bs{N} \ovl{X}_{\bs{N}} \ri ] \geq
\bb{ E } [ \bs{N} ] \mu$.  On the other hand, by the convexity of $- g(x)$, we have $\bb{ E } \li [ \bs{N} (- \ovl{X}_{\bs{N}} ) \ri ] \geq \bb{
E } [ \bs{N} ] (- \mu)$ or equivalently, $\bb{ E } \li [ \bs{N} \ovl{X}_{\bs{N}} \ri ] \leq \bb{ E } [ \bs{N} ] \mu$.  Hence, it must be true
that $\bb{ E } \li [ S_{\bs{N}} \ri ] = \bb{ E } \li [ \bs{N}  \ovl{X}_{\bs{N}}  \ri ] =  \bb{ E } [ \bs{N} ] \mu$, which is Wald's first
equation. Similarly, we can demonstrate that the inequality in the second assertion of Theorem \ref{genWald} is a generalization of Wald's
second equation.  As an illustration of the applications of Theorem \ref{genWald}, consider
\[
\bs{N} = \inf \li \{ n \in \mscr{N}: n > \f{1}{ g ( \ovl{X}_n ) },
\; g ( \ovl{X}_n ) > 0 \ri \}.
\]
Clearly, $\bs{N} > \f{1}{ g ( \ovl{X}_{\bs{N}} ) }, \; g (
\ovl{X}_{\bs{N}} ) > 0$  and thus  $\bs{N} g ( \ovl{X}_{\bs{N}} )
\geq 1$ almost surely provided that $\bb{E} [ \bs{N} ] < \iy$.  By
using the generalization of  Wald's first equation, we have $\bb{E}
[ \bs{N} ] g (\mu) \geq 1$, which implies the following result.

\beT \la{UseWald} Assume that $g$ is a concave function on a convex set $\mscr{D}$ in $\bb{R}^d$ such that $\mscr{D}$ contains $\mu$ and the
range of $\ovl{X}_{\bs{N}}$ and that $g(\mu) > 0$. Then, $\bb{E} [ \bs{N} ] \geq \f{1}{ g (\mu)}$.

\eeT

\subsection{Generalization of Lorden's Inequality}

In order to obtain tight bounds for stopping times, we need to
generalize Lorden's inequality \cite{Lorden}. In this direction, we
have obtained the following result.

 \beT \la{Lordengen}
 Let $Z_1, Z_2, \cd$ be i.i.d. random variables having the same distribution as $Z$ such that $\bb{E} [ Z^2 ] < \iy$.
 Assume that $\bs{\lm}$ is a random variable independent of $Z_i$ for all $i \in \bb{N}$.
 Define {\small $\mscr{M}_{\bs{\lm}} = \inf \li \{ n \in \bb{N}: \sum_{i = 1}^n Z_i > \bs{\lm} \ri
\}$} and $R_{\bs{\lm}} = \sum_{i =1}^{\mscr{M}_{\bs{\lm}}}  Z_i -
\bs{\lm}$.  Then, $\bb{E} [ R_{\bs{\lm}} ] \leq \f{ \bb{E} [ (Z^+)^2
] }{ \bb{E} [ Z ] } \Pr \{ Z <  \bs{\lm}  \} + \bb{E} [ (Z -
\bs{\lm} )^+ ]$.

\eeT

See Appendix \ref{Lordengenapp} for a proof.

In the following, we have extended Lorden's inequality to the case
that the increment of time indexes is not a constant.

\beT

\la{Lordenspec} Let $Z_1, Z_2, \cd$ be i.i.d. positive random
variables having the same distribution as $Z$ such that $\bb{E} [
Z^2 ] < \iy$.
 Assume that $\bs{\lm}$ is a random variable independent of $Z_i$ for all $i \in \bb{N}$.
 Define $\mscr{M}_{\bs{\lm}} = \inf \{ n \in \mscr{N}: \sum_{i = 1}^n Z_i > \bs{\lm} \}$
 and $R_{\bs{\lm}} = \sum_{i =1}^{\mscr{M}_{\bs{\lm}}}  Z_i - \bs{\lm}$.
  Define $Y = Z_1 + \cd + Z_{N_1}$ and $K = \max \{ N_{\ell + 1} - N_\ell: \ell \in \bb{N} \}$.
  Then, \[ \bb{E} [ R_{\bs{\lm}} ] \leq
\li ( (K - 1) \bb{E} [ Z ]  + \f{ \bb{E} [ Z^2 ] }{ \bb{E} [ Z ] }
\ri ) \Pr \{  Y  <  \bs{\lm}  \} + \bb{E} [ (Y - \bs{\lm} )^+ ].
\]

\eeT

See Appendix \ref{Lordenspecapp} for a proof.

\section{Bounding Average Stopping Time via Convex Optimization}

In this section, we shall demonstrate that the general problem
of bounding average stopping times can be converted into problems of convex
minimization, which can be readily solved by modern optimization theory and algorithms.
In some particular cases, it is possible to obtain
explicit bounds for average stopping times.

Consider the stopping time defined by (\ref{defFPT}).   To study the
properties of FPT, we introduce the concept of {\it deterministic
exit time} (DET).  The idea is to consider a deterministic motion
with two components, the first component is a one-dimensional motion
with constant velocity $1$, the other component is a $d$-dimensional
motion with constant velocity $v \in \bb{R}^d$. The overall motion
starts  from $(0, \bs{0}_d) \in \mscr{R}$. Clearly, the displacement
is $(t, v t)$ for $t \geq 0$. Depending on the structure of
$\mscr{R}$, the tuple $(t, v t)$ may or may not be contained in
$\mscr{R}$.  A positive real number $\tau$ is said to be a
 deterministic exit time for velocity $v$ if $(\tau, v \tau) \in \mscr{R}$ and there
exists a number $\de > 0$ such that $(\tau + \ep, v (\tau + \ep))
\notin \mscr{R}$  for $0 < \ep < \de$.  Define \be
\la{DEFDET66338899}
 \mscr{A}  (v) = \inf \{ t \geq 0: (t,  v t)
\notin \mscr{R} \}, \qqu \mscr{B}  (v) = \sup \{ t \geq 0: (t, v t)
\in \mscr{R} \} \ee for $v \in \bb{R}^d$.  In view of
(\ref{DEFDET66338899}), we call $\mscr{A} (v)$ the {\it infimum of
deterministic exit time} (IDET), and  $\mscr{B} (v)$  the {\it
suprimum of deterministic exit time} (SDET) for velocity $v$.
Clearly, if $\mscr{R}$ is a star domain, then the IDET and SDET are
equal, i.e., $\mscr{A} (v) = \mscr{B} (v)$. Particularly, this is
true when $\mscr{R}$ is convex.

\subsection{Lower Bound on Average Stopping Time}

   Regarding the lower bound of the expected value of the stopping time defined by (\ref{defFPT}),
   we have the following results.

\beT

\la{Convexlower}

Suppose that the stopping region $\mscr{R}^c$ is a convex set. Then,
$\bb{E} [ \bs{N} ] \geq \mscr{A} (\mu)$ provided that $\mscr{A}
(\mu) < \iy$.  Moreover, $\bb{E} [ \bs{N} ] = \iy$ provided that
$\mscr{A} (\mu) = \iy$. \eeT

See Appendix \ref{Convexlowerapp} for a proof.

\subsection{Upper Bounds on Average Stopping Time}

In order to develop upper bounds for the average of stopping time
$\bs{N}$ defined by (\ref{defFPT}), we shall investigate conditions
under which $\bb{E} [ \bs{N} ]$ is finite.  For this purpose,
consider the following conditions:

\bed

\item [(I)]  $\lim_{\ell \to \iy} \; \f{ N_{\ell + 1} }{ N_\ell } = 1$.

\item [(II)] $\mscr{R}$ is a closed convex set containing $(0, \bs{0}_d)$.

\item [(III)] $\mscr{B} (\mu) < \iy$.

\item [(IV)] Each element of $\bb{E} [ | X | ]$ is finite.

\eed

For the stopping time $\bs{N}$ defined by (\ref{defFPT}), we have
established the following result.

 \beT \la{BoundGen}

  If conditions (I)--(IV) are fulfilled, then $\bb{E} [\bs{N} ] < \iy$.

\eeT

See Appendix \ref{BoundGenapp} for a proof.    For the purpose of
bounding $\bb{E}[ \bs{N} ]$, we introduce the following conditions:

\bed

\item [(V)] There exist numbers $\lm > 0$ and $K$ such that $N_{\ell + 1} \leq \lm N_\ell + K$ for all $\ell \geq
0$.

\item [(VI)] $\{  (N_0, S_{N_0} )
\in \mscr{R} \}$ is a sure event.

\eed

It should be noted sample sizes used in group sequential methods
\cite{Bartroff, Jennsion, Mukhopadhyay, Proschan, Whitehead}
typically satisfy condition (V).   Define
\[
 \bs{M} = \sup \{ N_\ell: \ell \in \bb{Z}^+ , \; \bs{N} > N_\ell \}.
\]
Let $\bs{\ell}$ be the index of $N_\ell$ at the termination of the
random walk.  Then, $\bs{\ell}$ is a random variable such that
$\bs{N} = N_{\bs{\ell}}$ and $\bs{M} = N_{\bs{\ell} - 1}$.  If
$\bb{E} [\bs{N} ] < \iy$,  then it must be true that $\Pr \{ \bs{N}
< \iy \} = 1$.   It should be noted that $\bs{M}$ is not a stopping
time and thus $\bb{E} [ S_{\bs{M}} ]$ is, in general,  not equal to
$\bb{E} [ \bs{M} ] \mu$. In other words, Wald's first equation
\cite{Wald} is not applicable to $S_{\bs{M}}$, although it holds for
$S_{\bs{N}}$.   Clearly, as a consequence of the definition of
$\bs{M}$ and assumption (V), we have $\bs{N} \leq \lm \bs{M} + K$
and \be \la{indirect} \bb{E} [ \bs{N} ] \leq \lm \bb{E} [ \bs{M} ] +
K. \ee In view of (\ref{indirect}), to bound $\bb{E} [ \bs{N} ]$, it
suffices to bound $\bb{E} [ \bs{M} ]$.  For simplicity of notations,
define central moment
\[
\xi = \bb{E} [ | X - \mu | ],
\]
which will be used for obtaining upper bounds for $\bb{E} [ \bs{N}
]$. We have the following general result.

\beT

\la{GenConvex}  If conditions (I)--(VI) are fulfilled,  then  \be
\la{upper}
 \bb{E} [  \bs{M} ]
\leq \max_{ (t, s) \in \mscr{D}  }  t, \ee where $\mscr{D} = \li  \{
(t, s) \in \mscr{R} : | s - t \mu | \leq \f{1}{2} (\lm  t - r) \xi
\ri \}$ with $r = N_0 - K$.  \eeT

See Appendix \ref{GenConvexapp} for a proof.  It can be checked that $\mscr{D}$ is a convex set.  Moreover, $\max_{ (t, s) \in \mscr{D}  }  t =
- \min_{ (t, s) \in \mscr{D} } f(t,s)$, where $f(t,s) = - t$ is a convex function of $(t, s)$ contained in the convex set $\mscr{D}$.
Therefore, the upper bound in (\ref{upper}) can be readily evaluated by convex minimization.  With recent improvements in computing and in
optimization theory, convex minimization is nearly as straightforward as linear programming (see, e.g., \cite{Boyd} for a comprehensive
treatment).  Convex minimization problems can be solved by contemporary methods such as subgradient projection methods \cite{Polyak},
interior-point methods \cite{{Nesterov}}, etc.

In the case that $X$ is a bounded random vector, we have the following result.

\beT

\la{Convexsecond}  Suppose that $\Pr \{ \bs{a} \leq X \leq \bs{b} \}
= 1$, where $\bs{a}, \bs{b} \in \bb{R}^d$, and that conditions
(I)--(III), (V) and (VI) are fulfilled.  Define $\bs{v} = \f{ (\mu -
\bs{a}) (\bs{b} - \mu)  }{ \bs{b} - \bs{a} }$ and $r = N_0 - K$.
Then, \be \la{upperbound}
 \bb{E} [  \bs{M} ]
\leq \max_{ (t, s) \in \mscr{D}  }  t, \ee where $\mscr{D} = \{ (t,
s) \in \mscr{R} :  | s - t \mu | \leq (\lm t - r) \bs{v}, \; [ (\lm
- 1) t + K ] ( \mu - \bs{b} ) \leq s - t \mu \leq [ (\lm - 1) t + K
] ( \mu - \bs{a} ) \}$.  \eeT

See Appendix \ref{Convexsecondapp} for a proof.

In many situations, a stopping time is defined in terms of empirical
mean. Consider stopping time \be \la{STf}
 \bs{N} = \inf \{ n \in \mscr{N}: n >
g (\ovl{X}_n) \}. \ee For such stopping time, we have the following result.

\beT

\la{ConvexThird} Assume that $g$ is a concave function on a convex
set $D$ in $\bb{R}^d$ such that $\mu$ is an interior point of $D$
and that the range of $\ovl{X}_n$ is contained by $D$ for any $n \in
\{N_0, N_1, N_2, \cd \}$.  Assume that $K$ and $N_0$ are positive
integers such that $\{ N_0 \leq g(\ovl{X}_{N_0} ) \}$ is a sure
event and that $N_{\ell + 1} - N_\ell \leq K \leq N_0$ for $\ell \in
\bb{Z}^+$. Assume that each element of $\bb{E} [ | X | ]$ is finite.
Then,
\[
\bb{E} [ \bs{N} ] \leq K + \max_{ \se \in \mscr{D}  } g (\se),
\]
where $\mscr{D} = \li \{  \se \in D:  | \se - \mu | \leq \f{\xi}{2}
\ri \}$.

\eeT

See Appendix \ref{ConvexThirdapp} for a proof.

If $\mscr{N} = \bb{N}$, the stopping time defined by (\ref{STf})
becomes $\bs{N} = \inf \{ n \in \bb{N}: n >  g (\ovl{X}_n) \}$. For
such stopping time, we have the following result.

\beC

\la{Convexfour} Assume that $g$ is a non-negative concave function
on a convex set $D$ in $\bb{R}^d$ such that $\mu$ is an interior
point of $D$ and that the range of $\ovl{X}_n$ is contained by $D$
for any $n \in \bb{N}$.  Assume that each element of $\bb{E} [ | X |
]$ is finite. Then,
\[
\bb{E} [ \bs{N} ] \leq 2 + \max_{ \se \in \mscr{D}  } g (\se),
\]
where $\mscr{D} = \li \{  \se \in D:  | \se - \mu | \leq \f{\xi}{2}
\ri \}$.

\eeC

See Appendix \ref{Convexfourapp} for a proof.  In the case that $X$ is
a bounded random vector, we have the following result for the stopping
time defined by (\ref{STf}).

\beC

\la{try88} Assume that $\Pr \{ \bs{a} \leq X \leq \bs{b} \} = 1$,
where $\bs{a}, \bs{b} \in \bb{R}^d$. Assume that $g$ is a concave
function on a convex set $D$ in $\bb{R}^d$ such  that $\mu$ is an
interior point of $D$ and that the range of $\ovl{X}_n$ is contained
by $D$ for any $n \in \{N_0, N_1, N_2, \cd \}$. Assume that $K$ and
$N_0$ are positive integers such that  $\{ N_0 \leq g (\ovl{X}_{N_0}
) \}$ is a sure event and that $N_{\ell + 1} - N_\ell \leq K \leq
N_0$ for $\ell \in \bb{Z}^+$. Then,
\[
\bb{E} [ \bs{N} ] \leq K + \max_{ \se \in \mscr{D}  } g (\se),
\]
where $\mscr{D} = \{  \se \in D:  | \se - \mu | \leq \bs{v} \}$ with
$\bs{v} = \f{ (\mu - \bs{a}) (\bs{b} - \mu) }{ \bs{b} - \bs{a} }$.

\eeC

See Appendix \ref{try88app} for a proof.

\section{Bounding Average Stopping Time with Supporting Hyperplane}

In this section, we shall establish explicit bounds for the expected
value of the stopping time $\bs{N}$ defined by (\ref{defFPT}).  For
this purpose, we propose to use the concept of supporting hyperplane
to derive explicit bounds for $\bb{E} [ \bs{N} ]$.

Throughout this section,  we will assume that $\mscr{R}$ is a closed
convex set containing $(0, \bs{0}_d)$ and let $m = \mscr{B} (\mu)$
for the sake of notational simplicity.  As a consequence of the
convexity of $\mscr{R}$, there exists a supporting hyperplane of
$\mscr{R}$, passing through $(m, m \mu)$. Under the assumption that
the supporting hyperplane does not contain $(0, \bs{0}_d)$, such
supporting hyperplane can be expressed by the equation $A s + B t =
m$ with $A^\top \in \bb{R}^d$ and $B \in \bb{R}$. The supporting
hyperplane consists of points $(t, s)$ with $t \in \bb{R}$ and $s
\in \bb{R}^d$ satisfying the equation.    Specially, when the DET
function $\mscr{B} (v)$ is differentiable at $v = \mu$, the
supporting hyperplane is actually the tangent plane, which can be
explicitly expressed  in terms of the gradient of $\ln \mscr{B} (v)$
at $v = \mu$.   More specifically, let $\nabla(v)$ denote the
gradient of $\ln \mscr{B} (v)$ with respect to $v$, that is, \be
\la{defnabla}
 \nabla(v) = \f{\pa \ln \mscr{B} (v) } {\pa v} = \f{ 1 }{\mscr{B} (v)} \f{\pa \mscr{B} (v)   }{ \pa v }. \ee
By virtue of Lemma \ref{chenconvex}  in Appendix
\ref{discrereboundap}, the supporting hyperplane can be expressed as
\be \la{tangentplane}
 - V s + (1 + V \mu )t = m,
\ee where $V = \nabla(\mu)$ is the gradient of $\ln \mscr{B} (v)$ at
$v = \mu$.

To bound the expected value of the stopping time $\bs{N}$ defined by
(\ref{defFPT}), we have the following result.

\beT \la{lineartheorem}

Assume that conditions (I)--(VI) are fulfilled.  Let $A s + B t =
m$, where $m = \mscr{B} (\mu) > 0$, be a supporting hyperplane of
$\mscr{R}$ passing through $(m, m \mu)$. Assume that $\lm | A | \xi
< 2$. Then, $\bb{E} [ \bs{N} ] \leq \f{ \lm m + K }{ 1 - \f{1}{2}
\lm | A | \xi }$.  Moreover,  $\bb{E} [ \bs{N} ] \leq N_0 + \f{ \lm
m + K - N_0 }{ 1 - \f{1}{2} \lm | A | \xi }$ provided that $\lm m +
K \geq N_0$. \eeT

See Appendix \ref{lineartheoremapp} for a proof.  Clearly, if the
elements of $\xi$ are close to $0$,  then $\f{1}{2}  \lm | A | \xi$
is close to $0$ and the upper bound of $\bb{E} [ \bs{N} ]$ is close
to $\lm m + K$.

When the gradient of $\ln \mscr{B} (v)$ is available at $v = \mu$,
the following result can be readily derived from Theorem
\ref{lineartheorem}.

\beC

Assume that conditions (I)--(VI) are fulfilled.  Assume that
$\mscr{B} (v)$ is differentiable at $v = \mu$.  Let $V$ be the
gradient of $\ln \mscr{B} (v)$ at $v = \mu$.  Assume that $\lm | V |
\xi < 2$. Then, $\bb{E} [ \bs{N} ] \leq \f{ \lm m + K }{ 1 -
\f{1}{2} \lm | V | \xi }$.  Moreover,  $\bb{E} [ \bs{N} ] \leq N_0 +
\f{ \lm m + K - N_0 }{ 1 - \f{1}{2} \lm | V | \xi }$ provided that
$\lm m + K \geq N_0$.

\eeC

In situations that $X$ is a bounded random vector, we have applied
Theorem \ref{lineartheorem} to  obtain explicit bounds for $\bb{E} [
\bs{N} ] $ as follows.

\beC \la{lineartheorembounded}

Assume that $\Pr \{ \bs{a} \leq X \leq \bs{b} \} = 1$, where
$\bs{a}, \bs{b} \in \bb{R}^d$, and that conditions (I)--(III), (V)
and (VI) are fulfilled.  Let $A s + B t = m$, where $m = \mscr{B}
(\mu) > 0$, be a supporting hyperplane of $\mscr{R}$ passing through
$(m, m \mu)$. Define $\bs{v} = \f{ (\mu - \bs{a}) (\bs{b} - \mu) }{
\bs{b} - \bs{a} }$. Assume that $\lm | A | \bs{v} < 1$. Then,
$\bb{E} [ \bs{N} ] \leq \f{ \lm m + K }{ 1 - \lm | A | \bs{v} }$.
Moreover, $\bb{E} [ \bs{N} ] \leq N_0 + \f{ \lm m + K - N_0 }{ 1 -
\lm | A | \bs{v} }$ provided that $\lm m + K \geq N_0$.
 \eeC

To show Corollary \ref{lineartheorembounded}, note that $\xi \leq 2
\bs{v}$, as a result of Lemma \ref{convexupper} in Appendix
\ref{Convexsecondapp}.  Using this fact and  Theorem
\ref{lineartheorem} proves the corollary.

When the gradient of $\ln \mscr{B} (v)$ is available at $v = \mu$,
the following result is a direct consequence of  Corollary
\ref{lineartheorembounded}.

\beC

Assume that $\Pr \{ \bs{a} \leq X \leq \bs{b} \} = 1$, where
$\bs{a}, \bs{b} \in \bb{R}^d$, and that conditions (I)--(III), (V)
and (VI) are fulfilled.  Define $\bs{v} = \f{ (\mu - \bs{a}) (\bs{b}
- \mu) }{ \bs{b} - \bs{a} }$.  Assume that $\mscr{B} (v)$ is
differentiable at $v = \mu$.   Let $V$ be the gradient of $\ln
\mscr{B} (v)$ at $v = \mu$.  Assume that  $\lm | V | \bs{v} < 1$.
Then, $\bb{E} [ \bs{N} ] \leq \f{ \lm m + K }{ 1 - \lm | V | \bs{v}
}$. Moreover, $\bb{E} [ \bs{N} ] \leq N_0 + \f{ \lm m + K - N_0 }{ 1
- \lm | V | \bs{v} }$ provided that $\lm m + K \geq N_0$.

\eeC

In the sequel, we shall apply the concept of supporting hyperplane
and Lorden's inequality on overshoot to obtain explicit bounds for
average stopping times. Consider stopping time \be \la{newST8896}
\bs{N} = \li \{ n \in \bb{N}: \f{n}{K} \in \bb{N}, \;  (n, S_n)
\notin \mscr{R} \ri \}, \ee where  $K$ is a positive integer.  For
such stopping time, we have the following results.

\beT

\la{ChenLorden}

Assume that $\mscr{R}$ is a closed convex set containing $(0,
\bs{0}_d)$. Assume that each element of $\bb{E} [ | X |^2 ]$ is
finite.  Assume that there exists a supporting hyperplane $A s + B t
= m$, where $m = \mscr{B} (\mu) > 0$,  of $\mscr{R}$ passing through
$(m, m \mu)$.  The following assertions hold true.

(I):  $\bb{E}[ \bs{N} ]  \leq  m + K +  \bb{E} \li [ \li | A ( X -
\mu) \ri |^2 \ri ] \leq m + K +   || A ||^2 \times \bb{E} \li [ || X
- \mu ||^2 \ri ]$.

(II): If  the  elements of $X$ are mutually independent, then
$\bb{E}[ \bs{N} ]  \leq m + K +  A^2 \; \bb{E} [ (X - \mu)^2 ]$.

(III): If $\Pr \{ \bs{a} \leq X \leq \bs{b} \} = 1$, where $\bs{a},
\bs{b} \in \bb{R}^d$, then $\bb{E} [ \bs{N} ] \leq m + K ( u + v - u
v)$,  where $u = \f{1}{2} \li [ A (\bs{a} + \bs{b}) + |A| (\bs{a} -
\bs{b}) \ri ] + B$ and $v = \f{1}{2} \li [ A (\bs{a} + \bs{b}) + |A|
(\bs{b} - \bs{a}) \ri ] + B$.  In particular, $\bb{E} [ \bs{N} ]
\leq  m + K v^2  \li ( \f{ 1 - u }{ v - u } \ri )$ for $u < 0$.

 \eeT

See Appendix \ref{ChenLordenapp} for a proof.

When the gradient of $\ln \mscr{B} (v)$ is available at $v = \mu$,
we can apply  Theorem \ref{ChenLorden} to derive the following
explicit bounds for $\bb{E} [ \bs{N} ]$, where $\bs{N}$ is defined
by (\ref{newST8896}).

\beT \la{discrerebound} Assume that $\mscr{R}$ is a convex set
containing $(0, \bs{0}_d)$. Assume that $\bb{E} [||X||^2]$ is
finite. Assume that $\mscr{B} (v)$ is differentiable at $v = \mu$.
Let $V$ be the gradient of $\ln \mscr{B} (v)$ at $v = \mu$. Then,
\[
\bb{E} [ \bs{N} ] \leq \mscr{B} (\mu) + K + \bb{E} \li [  \li | V (
X - \mu) \ri |^2 \ri ] \leq \mscr{B} (\mu)  + K +   \li | \li | V
\ri | \ri |^2 \times
 \bb{E} \li [ \li | \li | X - \mu \ri | \ri |^2  \ri ].
\]
\eeT

See Appendix \ref{discrereboundap} for a proof.

We can apply Theorem \ref{discrerebound} to derive a simple bound
for the expectation of the first passage time for a random walk with
concave boundary.  More specifically, consider stopping time \be
\la{exam8888}
 \bs{N} = \inf \li \{ n \in \bb{N}: \f{n}{K} \in \bb{N}, \; S_n > f (n) \ri \},
\ee where $K$ is a positive integer and $S_n = \sum_{i = 1}^n X_i$
is the partial sum of i.i.d scalar random variables $X_1, X_2, \cd$,
which have the same distribution as $X$ with mean $\mu = \bb{E} [ X
]$ and variance $\si^2 = \bb{E} [ | X - \mu |^2 ] < \iy$.  We have
the following result.

\beC

\la{concavefunction} Assume that $f(t)$ is a concave function of $t
\in \bb{R}^+$ such that $f(0) > 0$. Assume that there exists a
positive number $m$ such that $m \mu = f (m)$.  Assume that $f(t)$
is differentiable at $t = m$.  Then, \be \la{vipformula} \bb{E} [
\bs{N} ] \leq m + K + \f{ \si^2  }{  | f^\prime (m) - \mu |^2 }, \ee
where $f^\prime (m)$ is the derivative of $f(t)$ at $t = m$.

\eeC

See Appendix \ref{concavefunctionapp} for a proof.

To use formula (\ref{vipformula}), we need to obtain $m$ from
equation $m \mu = f (m)$. In many cases, it is possible to derive an
explicit expression of $m$ from such equation. Even if $m$ cannot be
obtained analytically, it can still be readily computed by numerical
methods such as the bisection search method. Due to the concavity of
$f(.)$ and the existence of $m$ satisfying $m \mu = f (m)$, it must
be true that $t \mu > f(t)$ for large enough $t
> 0$. For example, we can find such value of $t$ as $2^k$ for some
integer $k
> 0$. Then, the number $m$ can be obtained by a bisection search
from interval $(0, 2^k)$.

\section{Bounding Average Stopping Time with Concentration
Inequalities} \la{boundconcentration88936}

In this section, we shall propose a method for bounding average
stopping times by virtue of concentration inequalities. Consider the
stopping time $\bs{N}$ defined by (\ref{defFPT}). Define  \[
 \de (t) = \inf \{ || v - \mu ||:  v \in \bb{R}^d, \; (t, v t)
 \in \mscr{R} \}
 \]
  for $t > 0$.  Define  $m = \mscr{B} (\mu)$ and
  \[
  \jmath = \min \{ i \in
\bb{N}:   N_i > m \}, \qqu \ka = \sup \{ i \in \bb{N}:  \Pr \{
\bs{N} > N_i \} > 0 \}.
\]
We have the following general result.

\beT \la{BoundConcentrationa88} $\bb{E} [ \bs{N} ] \leq N_1 +
\sum_{\ell = 1}^\ka  (N_{\ell + 1} - N_\ell ) \Pr \li \{ \li | \li |
\ovl{X}_{N_\ell} -  \mu \ri | \ri | \geq \de (N_\ell) \ri \}$.  \eeT

See Appendix \ref{BoundConcentrationa88app} for a proof.

 It should be noted that if $\mscr{R}$ is convex, then $\de (t)$
 can be readily obtained by convex minimization.  Making use of the concept
 of supporting hyperplane, we have the following result.

\beT  \la{BoundConcentrationb88}  Assume that the continuity region
$\mscr{R}$ is a closed convex set containing $(0, \bs{0}_d)$. Assume
that there exists a supporting hyperplane $A s + B t = m$, where $m
= \mscr{B} (\mu) > 0$,  of $\mscr{R}$ passing through $(m, m \mu)$.
Then, \be \la{inelinear8899883}
 \bb{E} [ \bs{N} ] \leq N_\jmath + \sum_{\ell =
\jmath  }^\ka  (N_{\ell + 1} - N_\ell )
 \Pr \li \{  || A || \times \li | \li | \ovl{X}_{N_\ell} -  \mu \ri | \ri | \geq 1 - \f{m}{N_\ell}  \ri \}. \ee

\eeT

See Appendix \ref{BoundConcentrationb88app} for a proof.

When the gradient of $\ln \mscr{B}(v)$ is available at $v = \mu$,
the supporting hyperplane is actually the tangent plane described by
(\ref{tangentplane}). Hence, applying Theorem
\ref{BoundConcentrationb88},   we have obtained the following
result.

\beC

Assume that the continuity region $\mscr{R}$ is a closed convex set
containing $(0, \bs{0}_d)$.  Assume that $\mscr{B} ( v )$ is
differentiable at $v = \mu$.  Let $V$ be the gradient of $\ln
\mscr{B} ( v )$  at $v = \mu$.  Then,
\[
\bb{E} [ \bs{N} ] \leq N_\jmath + \sum_{\ell = \jmath  }^\ka
(N_{\ell + 1} - N_\ell )
 \Pr \li \{  || V || \times \li | \li | \ovl{X}_{N_\ell} -  \mu \ri | \ri | \geq 1 - \f{m}{N_\ell}  \ri
 \}.
 \]

\eeC

In the case that $X$ is a scalar random  variable, it suffices to
use one-sided probabilistic inequality.

\beT

\la{scalarAFPT}

Assume that $X$ is a scalar random  variable.  Assume that the
continuity region $\mscr{R}$ is a closed convex set containing $(0,
0)$. Then, \bel &  & \bb{E} [ \bs{N} ] \leq N_\jmath + \sum_{\ell =
\jmath}^\ka (N_{\ell + 1} - N_\ell ) \Pr \li \{ \ovl{X}_{N_\ell}
\geq \mu + \de (N_\ell) \ri \} \qu \tx{if} \qu  \{
(t, s) \in \mscr{R}: t > m, \; s > \mu t   \} \neq \emptyset, \qqu \qqu \la{ineq2879a}\\
&  & \bb{E} [ \bs{N} ] \leq N_\jmath + \sum_{\ell = \jmath}^\ka
(N_{\ell + 1} - N_\ell ) \Pr \li \{ \ovl{X}_{N_\ell} \leq \mu - \de
(N_\ell) \ri \} \qu \tx{if} \qu  \{ (t, s) \in \mscr{R}: t > m, \; s
< \mu t \} \neq \emptyset. \qqu \qqu \la{ineq2879b} \eel Moreover,
\be \la{ineq2882}
 \bb{E} [ \bs{N} ] \leq N_\jmath + \sum_{\ell = \jmath
}^\ka (N_{\ell + 1} - N_\ell ) \Pr \li \{  A \li ( \ovl{X}_{N_\ell}
-  \mu \ri ) \geq  1 - \f{m}{N_\ell}   \ri \} \ee holds under
additional assumption that there exists a supporting hyperplane $A s
+ B t = m$ of $\mscr{R}$, where $m = \mscr{B} (\mu) > 0$,  passing
through $(m, m \mu)$.

\eeT

See Appendix \ref{scalarAFPTapp} for a proof.

If the gradient of $\ln \mscr{B}(v)$ is available at $v = \mu$, the
supporting hyperplane is actually the tangent plane described by
(\ref{tangentplane}). Hence, applying Theorem \ref{scalarAFPT},  we
have the following result.

\beC

\la{scalarAFPTV} Assume that $X$ is a scalar random  variable.
Assume that the continuity region $\mscr{R}$ is a closed convex set
containing $(0, 0)$.  Assume that $\mscr{B} ( v )$ is differentiable
at $v = \mu$.  Let $V$ be the gradient of $\ln \mscr{B} ( v )$  at
$v = \mu$.  Then, \[ \bb{E} [ \bs{N} ] \leq N_\jmath + \sum_{\ell =
\jmath }^\ka (N_{\ell + 1} - N_\ell )
 \Pr \li \{  V \li ( \ovl{X}_{N_\ell} -  \mu \ri )  \leq  \f{m}{N_\ell} - 1  \ri \}.
\]

\eeC

It should be noted that the probabilistic terms in Theorem
\ref{scalarAFPT} and Corollary \ref{scalarAFPTV} may be bounded by
concentration inequalities such as Chernoff bounds and Hoeffding
inequalities \cite{Chernoff, Hoeffding}.

\section{Bounds for Average Stopping Times of L\'{e}vy Processes}

In the preceding sections, our techniques for bounding stopping
times are devoted to discrete-time stochastic processes. Actually,
the principle of such techniques can be extended to continuous-time
stochastic processes.  L\'{e}vy processes is an important category
of stochastic processes in continuous time (see, \cite{Bertoin,
Sato} and the references therein).  In this section, we shall focus
on the problem of bounding stopping times pertaining to L\'{e}vy
processes.

Let $\{ X_t, \; t \geq 0 \}$ be a L\'{e}vy process on $\bb{R}^d$
such that $\bb{E} [X_1] = \mu$.  Making use of Theorem \ref{genJen},
we have obtained the following result.

\beT

\la{use66889933}

 Assume that $\bs{T}$ is a positive random variable such that $\bb{E} [ \bs{T} ] < \iy$ and that for
any possible value $t$ of $\bs{T}$, the event $\{ \bs{T} = t \}$
depends only on $\{X_s: 0 \leq s \leq t \}$.  Assume that $g$ is a
convex function on $\bb{R}^d$.  Then, $\bb{ E } \li [ \bs{T} g (
\ovl{X}_{\bs{T}} ) \ri ] \geq  \bb{ E } [ \bs{T} ] g (\mu)$.

\eeT

The proof of Theorem \ref{use66889933} is similar to that of Theorem
\ref{genWald}, which is given in Appendix \ref{genWaldapp}.   In the
sequel, we shall investigate stopping times associated with  the
L\'{e}vy process $\{ X_t, \; t \geq 0 \}$  by virtue of the
geometric convexity of the continuity or stopping region.  Define
 \[
X = X_1, \qqu \ovl{X}_t  = \f{ X_t }{t} \qqu \tx{for $t > 0$},
 \]
and \be \la {defFPTLevy} T = \inf \{ t > 0: (t, X_t) \notin \mscr{R}
\}. \ee In probability theory, $T$ is also called the first passage
time (FPT) for the continuous-time random walk.   As before, the
continuity region $\mscr{R}$ is a closed subset of $\{ (t, s): t \in
\bb{R}^+, \; s \in \bb{R}^d \}$ which contains $(0, \bs{0}_d)$.  The
complement of $\mscr{R}$, denoted by $\mscr{R}^c$,  is called the
stopping region.  Let $\mscr{A} (v)$ and $\mscr{B} (v)$ be IDET and
SDET functions defined by (\ref{DEFDET66338899}).

Regarding the finiteness of the expected value of $T$ defined by
(\ref{defFPTLevy}), we have the following result.

\beT

\la{them8899Levy}

 Assume that  $\mscr{R}$ is a convex set such that $\mscr{B}(\mu) < \iy$.
Then, $\bb{E} [ T ] < \iy$. \eeT

To prove Theorem \ref{them8899Levy}, consider stopping time $\bs{N}
= \inf \li \{ n \in  \bb{N}: (n, X_n)  \notin \mscr{R} \ri \}$.
Clearly, $T \leq \bs{N}$.  From Theorem \ref{BoundGen}, we know that
$\bb{E} [ \bs{N} ] < \iy$. Hence, $\bb{E} [ T ]  \leq \bb{E} [
\bs{N} ] < \iy$.

Regarding the lower bound of the expected value of the stopping time
$T$ defined by (\ref{defFPTLevy}), we have obtained  the following
results.

\beT  \la{Brown2}  Assume that the stopping region $\mscr{R}^c$ is a
convex set.  Then,  $\bb{E} [ T ] \geq \mscr{A} (\mu)$ provided that
$\mscr{A} (\mu) < \iy$.  Moreover,  $\bb{E} [ T ] = \iy$ provided
that $\mscr{A} (\mu) = \iy$.  \eeT

See Appendix \ref{Brown2app} for a proof.

  For stopping time $T = \inf \li \{ t > 0: t > \f{1}{g \li ( \ovl{X}_t \ri )}, \; g \li ( \ovl{X}_t \ri ) > 0 \ri \}$,
   we have derived the following result.

\beC \la{Brown4} Assume that $g$ is a concave function on $\bb{R}^d$
with $g(\mu) > 0$. Then, $\bb{E} [ T ] \geq \f{1}{g (\mu)}$. \eeC

See Appendix \ref{Brown4app} for a proof.

Making use of the concept of supporting hyperplane, we have obtained
upper bounds for the expected value of the stopping time $T$ defined
by (\ref{defFPTLevy})  as follows.

\beT \la{them8899l}

Assume that the second moment of each element of $X$ is finite.
Assume that $\mscr{R}$ is a convex set. Assume that there exists a
supporting hyperplane $A s + B t = \tau$, where $\tau = \mscr{B}
(\mu) > 0$,  of $\mscr{R}$ passing through $(\tau, \mu \tau)$. Then,
$\bb{E}[ T ]  \leq  \tau +  \bb{E} \li [ \li | A ( X - \mu) \ri |^2
\ri ] \leq \tau +   || A ||^2 \times \bb{E} \li [ || X - \mu ||^2
\ri ]$.

 \eeT

See Appendix \ref{them8899lpf} for a proof.  When the gradient of
$\ln \mscr{B}(v)$ at $v = \mu$ is $V$, if follows from Lemma
\ref{chenconvex} in Appendix \ref{discrereboundap} that  the
supporting hyperplane of $\mscr{R}$, passing through $(\tau, \mu
\tau)$ with $\tau = \mscr{B} (\mu) > 0$, is actually the tangent
plane $- V s + (1 + V \mu )t = \tau$.   Hence, applying Theorem
\ref{them8899l}, we have obtained the following results for bounding
the stopping time $T$ defined by (\ref{defFPTLevy}).

\beC

\la{ChenLordenLevySupport}

Assume that the second moment of each element of $X$ is finite.
Assume that $\mscr{R}$ is convex and that $\mscr{B} (v)$ is
differentiable at $v = \mu$.   Let $V$ be the gradient of $\ln
\mscr{B} (v)$ at $v = \mu$.  Then, $\bb{E} [ T ] \leq \tau + \bb{E}
\li [ | V ( X - \mu) |^2 \ri ] \leq \tau + \li | \li | V \ri | \ri
|^2 \times \bb{E} \li [ \li | \li | X - \mu \ri | \ri |^2 \ri ]$,
where $\tau = \mscr{B} (\mu)$.

\eeC

We can apply Corollary \ref{ChenLordenLevySupport} to derive a
simple bound for the expectation of the first passage time for a
L\'{e}vy process crossing a concave boundary. More specifically, let
$\{ X_t, \; t \geq 0 \}$ be a scalar L\'{e}vy process  such that
$\bb{E} [ X ] = \mu$ and $\bb{E} [ | X - \mu |^2 ] = \si^2$, where
$X = X_1$. Consider stopping time $T = \inf \li \{ t > 0:  X_t > f
(t) \ri \}$.  Making use of Corollary \ref{ChenLordenLevySupport}
and following a similar argument as that for Corollary
\ref{concavefunction} in Appendix \ref{concavefunctionapp}, we have
the following result.

\beC

\la{concavefunctionLevy} Assume that $f(t)$ is a concave function of
$t \in \bb{R}^+$ such that $f(0) > 0$. Assume that there exists a
positive number $\tau$ such that $\tau \mu = f (\tau)$.  Assume that
$f(t)$ is differentiable at $t = \tau$.  Then,  $\bb{E} [ T ] \leq
\tau + \f{ \si^2  }{  | f^\prime (\tau) - \mu |^2 }$,  where
$f^\prime (\tau)$ is the derivative of $f(t)$ at $t = \tau$.

\eeC

A special class of L\'{e}vy processes is the Brownian motion (see,
\cite{KARATZAS, Morters, Revuz} and the references therein).  Let
$\{ W_t, \; t \geq 0 \}$ be a Brownian motion on $\bb{R}^d$ with
drift coefficient $\mu$ such that $\bb{E} [ W_t ] = t \mu$ for $t
\geq 0$. Define $\ovl{W}_t  = \f{ W_t }{t}$ for $t
> 0$.  Consider stopping time $T = \inf \{ t > 0: (t, W_t) \notin
\mscr{R} \}$. We have the following result.

\beT \la{Brown1}

 Assume that $\mscr{R}$ is a convex set such that $\mscr{B} (\mu) < \iy$. Then, $\bb{E} [ T ]
\leq \mscr{B} (\mu) $.

\eeT

See Appendix \ref{Brown1app} for a proof.

For stopping time $T = \inf \{ t > 0: t > g \li ( \ovl{W}_t \ri )
\}$, we have the following result.

\beC

\la{Brown3}

Assume that $g$ is a concave function on $\bb{R}^d$ with $g(\mu) >
0$. Then, $\bb{E} [ T ] \leq g (\mu)$.  \eeC

See Appendix \ref{Brown3app} for a proof.

\bsk

Similar to Section \ref{boundconcentration88936}, we shall propose a
method for bounding average stopping times associated with L\'{e}vy
processes by virtue of concentration inequalities. Consider the
stopping time defined by (\ref{defFPTLevy}). Define  \[
 \de (t) = \inf \{ || v - \mu ||:  v \in \bb{R}^d, \; (t,  v t)
 \in \mscr{R} \} \qqu \tx{for $t > 0$},
 \]
\[
\tau = \mscr{B} (\mu), \qqu  c = \sup \{ t \in \bb{R}^+: \Pr \{ T
> c \} > 0 \}.
\]
For the purpose of bounding $\bb{E} [ T ]$, we use Legesgue
integration in all bounds for $\bb{E} [ T ]$ in the remainder of
this section.
 We have the following general result.

\beT \la{BoundConcentrationaLevya} $\bb{E} [ T ] \leq \int_{0}^c \Pr
\li \{  \li | \li | \ovl{X}_t - \mu \ri | \ri | \geq \de (t) \ri \}
dt$.  \eeT

See Appendix \ref{BoundConcentrationaLevyaapp} for a proof.

 It should be noted that for all $t \in \bb{R}^+$, if the continuity region $\mscr{R}$ is convex, then  $\de (t)$
 can be readily obtained by convex minimization.
 Mimicking the proof of Theorem \ref{BoundConcentrationb88}, we have established the following results.

\beT  \la{BoundConcentrationLevyb} Assume that  the continuity
region $\mscr{R}$ is a closed convex set containing $(0, \bs{0}_d)$.
Assume that there exists a supporting hyperplane $A s + B t = \tau$,
where $\tau = \mscr{B} (\mu) > 0$,  of $\mscr{R}$ passing through
$(\tau, \mu \tau)$. Then,  $\bb{E} [ T ] \leq \tau + \int_\tau^c \Pr
\li \{  || A || \times \li | \li | \ovl{X}_{t} -  \mu \ri | \ri |
\geq  1 - \f{\tau}{t} \ri \} dt$.   \eeT

When the gradient of $\ln \mscr{B}(v)$ at $v = \mu$ assumes value
$V$, the supporting hyperplane is actually the tangent plane $- V s
+ (1 + V \mu) t = \tau$. Hence, applying Theorem
\ref{BoundConcentrationLevyb}, we have obtained the following
result.

\beC

Assume that  the continuity region $\mscr{R}$ is a closed convex set
containing $(0, \bs{0}_d)$.  Assume that $\mscr{B} ( v )$ is
differentiable at $v = \mu$.  Let $V$ be the gradient of $\ln
\mscr{B} ( v )$ at $v = \mu$. Then,
\[
\bb{E} [ T ] \leq \tau + \int_\tau^c
 \Pr \li \{  || V || \times \li | \li | \ovl{X}_t -  \mu \ri | \ri | \geq 1 - \f{\tau}{t}  \ri
 \} dt.
 \]
 \eeC

Mimicking the proof of Theorem \ref{scalarAFPT}, we have shown the
following result.

\beT

\la{scalarAFPTLevy}

Assume that $X$ is a scalar random  variable. Assume that  the
continuity region $\mscr{R}$ is a closed convex set containing $(0,
0)$. Then, \bee &  & \bb{E} [ T ] \leq \tau + \int_\tau^c \Pr \li \{
\ovl{X}_t \geq \mu + \de (t) \ri \} dt \qu \tx{if} \qu \{ (t, s) \in
\mscr{R}: t > \tau, \; s
> \mu t   \} \neq
\emptyset,  \\
&  & \bb{E} [ T ] \leq \tau + \int_\tau^c \Pr \li \{ \ovl{X}_t \leq
\mu - \de (t) \ri \} dt \qu \tx{if} \qu \{ (t, s) \in \mscr{R}: t >
\tau, \; s < \mu t   \} \neq \emptyset. \eee Moreover, $\bb{E} [ T ]
\leq \tau + \int_\tau^c  \Pr \li \{  A \li ( \ovl{X}_t -  \mu \ri )
\geq  1 - \f{\tau}{t}   \ri \} dt$ holds under additional assumption
that there exists a supporting hyperplane $A s + B t = \tau$, where
$\tau = \mscr{B} (\mu) > 0$,  of $\mscr{R}$ passing through $(\tau,
\tau \mu)$.
 \eeT

If the gradient of $\ln \mscr{B}(v)$ at $v = \mu$ assumes value $V$,
the supporting hyperplane is actually the tangent plane $- V s + (1
+ V \mu) t = \tau$. Hence, applying Theorem \ref{scalarAFPTLevy}, we
have the following result.

\beC

Assume that $X$ is a scalar random  variable. Assume that  the
continuity region $\mscr{R}$ is a closed convex set containing $(0,
0)$.  Assume that $\mscr{B} ( v )$ is differentiable at $v = \mu$.
Let $V$ be the gradient of $\ln \mscr{B} ( v )$ at $v = \mu$. Then,
$\bb{E} [ T ] \leq \tau + \int_\tau^c \Pr \li \{  V \li ( \ovl{X}_t
- \mu \ri ) \leq \f{\tau}{t} - 1  \ri \} dt$.  \eeC

\section{Conclusion}

In this paper, we have established a geometric approach for bounding
average stopping times.   The central idea of our approach is to
explore the geometric convexity of the continuity or stopping
regions.  Our approach are effective for a wide variety of stopping
times which involve random vectors, nonlinear boundary, constraint
of time indexes, etc.   Tight bounds are obtained for stopping times
in a general setting, which are explicit or readily computable.  A
probabilistic characterization is established for convex sets.
Extensions are developed for classical results such as Jensen's
inequality, Wald's equations and Lorden's inequality.

\appendix

\section{Proof of Theorem \ref{fundamental} } \la{fundamentalapp}

We need some preliminary results.  If $X$ is a random variable such
that $\Pr \{ X < c \} = 1$, then it is intuitive that $\bb{E} [ X ]
< c$. However, there exists no proof in the literature for such
intuition. Since the strictness of the inequality plays a crucial
role in our proof of the theorem, we provide a rigorous proof in the
sequel.

\beL \la{simpleb} If $X$ is a scalar random variable such that $\Pr
\{ X < c \} = 1$, then $\bb{E} [ X ] < c$.  Similarly, if $X$ is a
scalar  random variable such that $\Pr \{ X > c \} = 1$,  then
$\bb{E} [ X ]
> c$. \eeL

\bpf

We claim that there exists a positive number $\vep > 0$ such that
$\Pr \{ X \leq c - \vep \} > 0$. To prove the claim, we use a
contradiction method. Suppose that the claim is not true.  Then,
$\Pr \{ X \leq c - \vep \} = 0$ for any $\vep> 0$. It follows that
 $\Pr \{ X < c \} = \lim_{ \vep \downarrow 0 } \Pr \{ X \leq c - \vep
\} = 0$.  This contradicts to the assumption that $\Pr \{ X < c \} =
1$.  So, we have proved the claim.

Now let $\vep > 0$ be a positive number such that $\Pr \{ X \leq c - \vep \} > 0$.
Since $\Pr \{ X < c \} = 1$, we have  \bee \bb{E} [ X ] & =
& \bb{E} \li [ X \; \bb{I}_{ \{
X \leq c - \vep \} } \ri ] + \bb{E} \li [ X \; \bb{I}_{\{ c - \vep < X < c \} } \ri ]\\
& \leq & (c - \vep) \Pr \{  X \leq c - \vep \} + c  \Pr \{ c - \vep < X < c  \} \\
& = & (c - \vep) \Pr \{  X \leq c - \vep \} + c ( 1 - \Pr \{  X \leq c - \vep \} )\\
& = &  - \vep \Pr \{  X \leq c - \vep \} + c < c.     \eee
This proves the first assertion.   The second assertion can be shown in a similar
way.

\epf

\beL

\la{closed889}

Assume that $D$ is a closed convex set and $\bs{\mcal{X}}$ is a random vector
such that $\Pr \{ \bs{\mcal{X}} \in D \} = 1$, then $\bb{E} [
\bs{\mcal{X}} ] \in D$. \eeL

\bpf

 We shall use a contradiction method. Denote $\mu = \bb{E} [ \bs{\mcal{X}} ] $.
 Suppose $\mu \notin D$, i.e., $\mu$ is an exterior point
of $D$.  By the hyperplane separation theorem \cite[Theorem 4.11,
page 170]{Bertsimas}, there exists a row vector $\bs{\al}$ such that
$\bs{\al} \mu  <  \bs{\al} Z$ for all $Z \in D$.  Since $\Pr \{
\bs{\mcal{X}} \in D \} = 1$, it must be true that $\Pr \{ \bs{\al}
\mu < \bs{\al} \bs{\mcal{X}} \} = 1$.  Hence,  $\Pr \{ \bs{\al} \mu
- \bs{\al} \bs{\mcal{X}} < 0 \} = 1$.  It follows from Lemma
\ref{simpleb} that $\bb{E} [ \bs{\al} \mu - \bs{\al} \bs{\mcal{X}} ]
< 0$, which implies that
\[
\bs{\al} \mu  < \bb{E} [ \bs{\al} \bs{\mcal{X}} ] = \bs{\al} \bb{E}
[ \bs{\mcal{X}} ] =  \bs{\al} \mu.
\]
This is a contradiction.  Therefore, it must be true that $\mu \in
D$. The proof of the lemma is thus completed.

\epf

\beL

\la{simplea} If $X$ is a scalar random variable such that $0 < \Pr
\{ X < 0 \} \leq \Pr \{ X \leq 0 \} = 1$.  Then, $\bb{E} [ X ] < 0$.

\eeL

\bpf

We claim that there exists a positive number $\vep > 0$ such that
$\Pr \{ X \leq - \vep \} > 0$. To prove the claim, we use a
contradiction method. Suppose that the claim is not true.  Then,
$\Pr \{ X \leq - \vep \} = 0$ for any $\vep> 0$.  It follows that
 $\Pr \{ X < 0 \} = \lim_{ \vep \downarrow 0 } \Pr \{ X \leq - \vep \}
= 0$.  This contradicts to the assumption that $\Pr \{ X < 0 \} >
0$. So, we have proved the claim.

Now let $\vep > 0$ be a positive number such that $\Pr \{ X \leq - \vep \} > 0$.
Since $\Pr \{ X \leq 0 \} = 1$, we have \bee \bb{E} [ X ] & =
& \bb{E} \li [ X \; \bb{I}_{\{
X \leq - \vep \} } \ri ] + \bb{E} \li [ X \; \bb{I}_{ \{ - \vep < X \leq 0 \} } \ri ]\\
& \leq & - \vep \Pr \{  X \leq - \vep \} + 0 \times \Pr \{  - \vep < X \leq 0  \} \\
& \leq & - \vep \Pr \{  X \leq - \vep \} < 0.  \eee This completes the proof of the lemma.

\epf

\beL

\la{closed8999}

Assume that $D$ is a closed convex set and $\bs{\mcal{X}}$ is a
random vector such that $\Pr \{ \bs{\mcal{X}} \in D \} = 1$ and $\mu
= \bb{E} [ \bs{\mcal{X}} ] \in \pa D$, then there exist a nonzero
row vector $\bs{\al}$ and a constant $\ba$ such that $\Pr \{
\bs{\al} \bs{\mcal{X}}  + \ba = 0 \} = 1$. \eeL

\bpf

As a consequence of the convexity of $D$ and the assumption that
$\mu = \bb{E} [ \bs{\mcal{X}} ] \in \pa D$, it is possible to
construct a supporting hyperplane $ \bs{\al} Z + \ba = 0$ through
$\mu$, where $\bs{\al}$ is a nonzero row vector and $\ba$ is a
constant, such that $\bs{\al} Z + \ba \leq 0$ for all $Z \in D$. By
the assumption that $\Pr \{ \bs{\mcal{X}} \in D \} = 1$, we have
\[
\Pr \{ \bs{\al} \bs{\mcal{X}}  + \ba \leq 0  \} = 1.
\]
Since $\mu$ is in the supporting hyperplane, we have $ \bs{\al}
\bb{E} [ \bs{\mcal{X}} ]  + \ba = 0$.  We claim that $\Pr \{
\bs{\al} \bs{\mcal{X}}  + \ba = 0 \} = 1$. To prove this claim, we
use a contradiction method.  Suppose the claim is not true. Then,
\[ 0 < \Pr \{ \bs{\al} \bs{\mcal{X}}  + \ba < 0 \} \leq \Pr \{
\bs{\al} \bs{\mcal{X}} + \ba \leq 0 \} = 1. \] It follows from Lemma
\ref{simplea} that $\bb{E} [ \bs{\al} \bs{\mcal{X}} + \ba ] < 0$.
This implies that $\bs{\al} \bb{E} [ \bs{\mcal{X}} ]  + \ba = \bb{E}
[ \bs{\al} \bs{\mcal{X}}  + \ba ] < 0$,  which contradicts to the
fact that $\bs{\al} \bb{E} [ \bs{\mcal{X}} ]  + \ba = 0$.  The claim
is thus established. Hence, it must be true that $\Pr \{
 \bs{\al} \bs{\mcal{X}}  + \ba = 0  \} = 1$.   This completes the proof of the lemma.

 \epf

We are now in a position to prove the theorem.  We shall argue by a
mathematical induction on the dimension $n$ of $\mscr{D}$. For the
dimension $n = 1$, the convex set $\mscr{D}$ must be an interval of
the form $\mscr{D} = [a, b]$, or $\mscr{D} = (a, b)$, or $\mscr{D} =
[a, b)$, or $\mscr{D} = (a, b]$. Making use of Lemma \ref{simpleb},
it is easy to see $\bb{E} [ \bs{\mcal{X}} ] \in \mscr{D}$ as a
consequence of $\Pr \{ \bs{\mcal{X}} \in \mscr{D} \} = 1$. Suppose
the conclusion $\bb{E} [ \bs{\mcal{X}} ] \in \mscr{D}$ holds for
dimension $n - 1$. To complete the induction process, we need to
show, based on such hypothesis,  that the inclusion relationship  $
\bb{E} [ \bs{\mcal{X}} ] \in \mscr{D}$ holds for dimension $n$.  Let
$\overline{\mscr{D}}$ denotes the closure of $\mscr{D}$. By Lemma
\ref{closed889}, we have shown $ \bb{E} [ \bs{\mcal{X}} ] \in
\ovl{\mscr{D}}$. If $\mu = \bb{E} [ \bs{\mcal{X}} ]$ is not
contained in the boundary of $\ovl{\mscr{D}}$, then it must be true
that $\mu \in \mscr{D}$. Hence, to show $ \bb{E} [ \bs{\mcal{X}} ]
\in \mscr{D}$ for dimension $n$, it suffices to show it under the
assumption that $\mu = \bb{E} [ \bs{\mcal{X}} ]$ is contained in the
boundary of $\ovl{\mscr{D}}$. We proceed as follows.  Making use of
Lemma \ref{closed8999} and the assumption that $\mu = \bb{E} [
\bs{\mcal{X}} ]$ is contained in the boundary of $\ovl{\mscr{D}}$,
we conclude that there exist a nonzero row vector $\bs{\al}$ and a
constant $\ba$ such that $\Pr \{ \bs{\al} \bs{\mcal{X}}  + \ba = 0
\} = 1$. Define
\[
\mscr{S} =  \{ Z \in \mscr{D}:  \bs{\al} Z  + \ba = 0 \}.
\]
Then, $\mscr{S}$ is convex and $\Pr \{  \bs{\mcal{X}} \in  \mscr{S}
\} = 1$.  Without loss of any generality, assume that the $i$-th
element of $\bs{\al}$, denoted by $\al_i$, is nonzero. Define a
linear transform $\mscr{T}: \mscr{S}  \mapsto D$  such that for
every element $Z = [ z_1, \cd, z_n ]^\top$ in $\mscr{S}$, there
exists a corresponding vector $U = [ u_1, \cd, u_n ]^\top =
\mscr{T}(Z)$ such that
\[
u_i =  \bs{\al} Z  + \ba, \qqu u_\ell = z_\ell, \qqu \ell \in \{ 1,
\cd, n \} \setminus \{i\}
\]
or equivalently,  \be \la{trans88}
 U =  (I +  \bs{e}_i \bs{\al} - \bs{e}_i \bs{e}_i^\top  ) Z + \ba \bs{e}_i,
\ee where $I$ is an identity matrix of size $n \times n$ and
$\bs{e}_i$ is a column matrix with all elements being $0$ except the
$i$-th element being $1$. Note that $D = \{ \mscr{T} (Z): Z \in
\mscr{S} \}$ must be convex because the transform $\mscr{T}$ is
linear and $\mscr{S}$ is convex.  Define $\bs{Y} = [\bs{y}_1, \cd,
\bs{y}_n ]^\top = \mscr{T} (\bs{\mcal{X}} )$. Then, \[ \Pr \{ \bs{Y}
\in D \} = 1, \qqu \Pr \{ \bs{y}_i = 0 \} = \Pr \{ \bs{\al}
\bs{\mcal{X}} + \ba = 0  \} = 1 \]
 and $\bb{E} [ \bs{y}_i ] = 0$. Define \[
 D^* = \{ [u_1, \cd, u_i, u_{i+1}, \cd, u_{n}]^\top: [u_1, \cd,
u_{n}]^\top \in D \}. \]
 Then, $D^*$ is convex because $D$ is convex.
 Define random vector $\bs{V} = [ \bs{v}_1, \cd, \bs{v}_{n-1} ]^\top$ such that
$\bs{v}_\ell = \bs{y}_\ell, \; \ell = 1, \cd, i - 1$ and
$\bs{v}_\ell = \bs{y}_{\ell + 1}, \; \ell = i, \cd, n - 1$. Then, $\Pr \{ \bs{V} \in
D^* \} = 1$.  Since $D^*$ is a convex set of $(n-1)$ dimension
and  $\Pr \{ \bs{V} \in D^* \} = 1$, it follows from the induction hypothesis
that $\bb{E} [ \bs{V} ] \in D^*$.   This implies that $\bb{E} [ \bs{Y} ] \in D$.

It can be checked that the determinant of the matrix $I +  \bs{e}_i
\bs{\al} -  \bs{e}_i \bs{e}_i^\top$ in (\ref{trans88}) is equal to
$\al_i$, which is  nonzero. Hence, $I +  \bs{e}_i \bs{\al} -
 \bs{e}_i \bs{e}_i^\top$ is invertible, and it follows that
\[
Z = (I +  \bs{e}_i \bs{\al} -  \bs{e}_i \bs{e}_i^\top )^{-1} (U -
\ba \bs{e}_i ).
\]
This implies that the transform $\mscr{T}$ is a one-to-one mapping
from $\mscr{S}$ to $D$ and thus the transform is invertible. Note
that $\bb{E} [\bs{Y}] = \mscr{T} ( \bb{E} [ \bs{\mcal{X}} ] )$ and
the transform $\mscr{T}$ maps $\mscr{S}$ into $D$. Now, we have
$\bb{E} [ \bs{Y} ] \in D$. Taking the inverse transform of
$\mscr{T}$ yields $\bb{E} [ \bs{\mcal{X}} ] \in \mscr{S} \subseteq
\mscr{D}$.   This completes the process of the mathematical
induction and the theorem is thus established.

\sect{Proof of Theorem \ref{genJen} } \la{genJenapp}

We need some preliminary result.

\beL
 \la{BoundaryCon}
Suppose that $g(z)$ is a  convex function of $z \in \mscr{D}$, where
$\mscr{D}$ is a convex set  in $\bb{R}^n$.
 Define $f (t, s) = t g \li (\f{s}{t} \ri )$ for $t \neq 0$ and $s$ such that $\f{s}{t} \in \mscr{D}$. Then, $f(t, s)$
 is a  convex function of $t > 0$ and $s$ such that $\f{s}{t} \in \mscr{D}$.  Similarly, $f(t, s)$
 is a  concave function of $t < 0$ and $s$ such that $\f{s}{t} \in \mscr{D}$.

\eeL

\bpf To show the first assertion, it suffices to show that the
inequality $f ( \sum_{\ell = 1}^k \lm_\ell t_\ell, \; \sum_{\ell =
1}^k \lm_\ell s_\ell  ) \leq \sum_{\ell = 1}^k \lm_\ell f (t_\ell,
s_\ell)$ holds for any $(t_\ell, s_\ell), \; \ell = 1, \cd, k$ such
that $t_\ell
> 0, \; \f{s_\ell}{t_\ell} \in \mscr{D}$ and nonnegative numbers
$\lm_\ell, \; \ell = 1, \cd, k$ such that $\sum_{\ell = 1}^k
\lm_\ell = 1$.    Define $A = \sum_{\ell = 1}^k \lm_\ell t_\ell$ and
$\ro_\ell = \f{\lm_\ell t_\ell}{ A}$ for $\ell = 1, \cd, k$.  Since
$\ro_\ell, \; \ell = 1, \cd, k$ are nonnegative numbers satisfying
$\sum_{\ell = 1}^k \ro_\ell = 1$ and the function $g$ is convex, we
have
\[
\sum_{\ell = 1}^k \ro_\ell \; g \li (\f{s_\ell}{t_\ell} \ri ) \geq
g \li ( \sum_{\ell = 1}^k \ro_\ell \; \f{s_\ell}{t_\ell} \ri ) = g
\li ( \sum_{\ell = 1}^k \f{\lm_\ell t_\ell}{ A} \;
\f{s_\ell}{t_\ell} \ri ) = g \li (  \f{ \sum_{\ell = 1}^k \lm_\ell
s_\ell}{ A} \;
 \ri ).
\]
It follows that \bee \sum_{\ell = 1}^k \lm_\ell  f(t_\ell, s_\ell)
& = & \sum_{\ell = 1}^k \lm_\ell t_\ell g \li (\f{s_\ell}{t_\ell} \ri )
 =  A \sum_{\ell = 1}^k \ro_\ell \; g \li (\f{s_\ell}{t_\ell} \ri )\\
& \geq &   A g \li (  \f{ \sum_{\ell = 1}^k \lm_\ell s_\ell}{ A} \;
 \ri ) =  \li ( \sum_{\ell = 1}^k \lm_\ell t_\ell \ri ) g \li
(\f{\sum_{\ell = 1}^k \lm_\ell s_\ell}{ \sum_{\ell = 1}^k \lm_\ell
t_\ell } \ri ) = f \li ( \sum_{\ell = 1}^k \lm_\ell t_\ell, \;
\sum_{\ell = 1}^k \lm_\ell s_\ell \ri ).  \eee This proves the first
assertion. The second assertion can be shown in a similar way.

\epf

We shall only show the first assertion, since the second assertion
can be shown in a similar way.  Define $f(t, s) = t g \li ( \f{s}{t}
\ri )$. Since $g(z)$ is a convex function of $z \in \mscr{D}$,  it
follows from Lemma \ref{BoundaryCon} that $f(t, s)$ is a convex
function of $t > 0$ and vector $s \in \bb{R}^n$ such that $\f{s}{t}
\in \mscr{D}$. Hence,  there exist a row vector $\bs{\al}$ and
number $\ba$ such that
\[
f (t, s) \geq f ( \bb{E} [ Y ], \bb{E} [ \bs{Z} ] ) +  \bs{\al} (s -
\bb{E} [ \bs{Z} ]) + \ba (t - \bb{E} [ Y ])
\]
for $t > 0$ and vector $s \in \bb{R}^n$ such that $\f{s}{t} \in
\mscr{D}$. As a consequence of this result and the assumption that
$Y > 0, \; \f{\bs{Z}}{Y} \in \mscr{D}, \; \f{ \bb{E} [ \bs{Z} ] }{
\bb{E} [ Y ] } \in  \mscr{D}$, we have
\[
f (Y, \bs{Z}) \geq f ( \bb{E} [ Y ], \bb{E} [ \bs{Z} ] ) + \bs{\al}
(\bs{Z} - \bb{E} [ \bs{Z} ])  + \ba (Y - \bb{E} [ Y ]).
\]
Applying the definition of the function $f$ to the above inequality yields
\[
Y g \li ( \f{\bs{Z}}{Y} \ri ) \geq \bb{E} [ Y ] g \li ( \f{ \bb{E} [
\bs{Z} ] }{ \bb{E} [ Y ] } \ri ) +  \bs{\al} (\bs{Z} - \bb{E} [
\bs{Z} ])  + \ba (Y - \bb{E} [ Y ]).
\]
Taking expectations on both sides leads to \[ \bb{E} \li [ Y g \li (
\f{\bs{Z}}{Y} \ri ) \ri ] \geq \bb{E} [ Y ] g \li ( \f{ \bb{E} [
\bs{Z} ] }{ \bb{E} [ Y ] } \ri ) +  \bs{\al} \bb{E} [ \bs{Z} -
\bb{E} [ \bs{Z} ] ]  + \ba \bb{E} [ Y - \bb{E} [ Y ] ] = \bb{E} [ Y
] g \li ( \f{ \bb{E} [ \bs{Z} ] }{ \bb{E} [ Y ] } \ri ).
\]

\sect{Proof of Theorem \ref{genWald} } \la{genWaldapp}

To show the first assertion, we can use the first inequality of
Theorem \ref{genJen} to conclude that \bee \bb{ E } \li [ \bs{N} g (
\ovl{X}_{\bs{N}} ) \ri ]  =  \bb{ E } \li [ \bs{N} g \li (
\f{1}{\bs{N}} \sum_{i=1}^{\bs{N}} X_i \ri ) \ri ] \geq \bb{E} [
\bs{N} ] g \li ( \f{1}{ \bb{E} [ \bs{N} ]} \bb{E} \li [
\sum_{i=1}^{\bs{N}} X_i \ri ] \ri ). \eee By virtue of Wald's first
equation, we have $\bb{E} \li [ \sum_{i=1}^{\bs{N}} X_i \ri ] =
\bb{E} [ \bs{N} ] \mu$.  Hence,
\[
\bb{ E } \li [ \bs{N} g ( \ovl{X}_{\bs{N}} ) \ri ]
\geq \bb{E} [  \bs{N} ] g \li ( \f{1}{ \bb{E} [ \bs{N} ]} \bb{E} [ \bs{N} ] \mu \ri ) =
\bb{E} [  \bs{N} ] g (\mu).
\]

To show the second assertion, we can use the first inequality of
Theorem \ref{genJen} to conclude that \bee \bb{ E } \li [ \bs{N} g (
\ovl{V}_{\bs{N}} ) \ri ]  =  \bb{ E } \li [ \bs{N} g \li (
\f{1}{\bs{N} } \li ( \sum_{i=1}^{\bs{N}} X_i - \bs{N} \mu \ri )^2
\ri ) \ri ] \geq \bb{E} [  \bs{N} ] g \li ( \f{1}{ \bb{E} [ \bs{N}
]} \bb{E} \li [ \li ( \sum_{i=1}^{\bs{N}} X_i - \bs{N} \mu \ri )^2
\ri ] \ri ). \eee By virtue of Wald's second equation, we have
$\bb{E} \li [ \li ( \sum_{i=1}^{\bs{N}} X_i - \bs{N} \mu \ri )^2 \ri
]  = \bb{E} [ \bs{N} ] \nu$.  Hence,
\[
\bb{ E } \li [ \bs{N} g ( \ovl{V}_{\bs{N}} ) \ri ]
\geq \bb{E} [  \bs{N} ] g \li ( \f{1}{ \bb{E} [ \bs{N} ]} \bb{E} [ \bs{N} ] \nu \ri ) =
\bb{E} [ \bs{N} ] g (\nu).
\]

\sect{Proof of Theorem \ref{Lordengen} } \la{Lordengenapp}

Define $\ze = \bs{\lm} - Z_1$.  Let $F_\ze (.)$ denotes the
cumulative distribution of $\ze$.  Note that \bel  \bb{E} \li [ \li
( \sum_{i =2}^{\mscr{M}_{\bs{\lm}}}  Z_i - ( \bs{\lm} - Z_1) \ri )
\; \bb{I}_{ \{ Z_1 < \bs{\lm} \} } \ri ]
&  =  & \bb{E} \li [ \li ( \sum_{i =2}^{\mscr{M}_{\bs{\lm}}}  Z_i - \ze \ri ) \; \bb{I}_{ \{ \ze > 0 \} } \ri ] \nonumber\\
& = & \int_{u > 0} \bb{E} \li [ \li ( \sum_{i
=2}^{\mscr{M}_{\bs{\lm}}}  Z_i - \ze \ri ) \; \mid \ze = u  \ri ] d
F_\ze (u).  \la{use3698a} \eel By the definition of
$\mscr{M}_{\bs{\lm}}$, we have \be \la{use3698b} \bb{E} \li [ \li (
\sum_{i =2}^{\mscr{M}_{\bs{\lm}}}  Z_i - \ze \ri ) \; \mid \ze = u
\ri ] = \bb{E} \li [ \li ( \sum_{i =2}^{\mcal{M}_{u}} Z_i - u \ri )
\; \mid \ze = u  \ri ], \ee where $\mcal{M}_u = \inf \li \{ n \geq
2: \sum_{i=2}^n Z_i  > u \ri  \}$.  Since the random variables $Z_1,
Z_2, \cd$ and $\bs{\lm}$ are independent, it follows that $\ze$ and
$Z_2, Z_3, \cd$ are independent.  Hence,  \be \la{use3698c}
 \bb{E} \li [ \li ( \sum_{i =2}^{\mcal{M}_{u}} Z_i - u \ri ) \; \mid \ze = u  \ri ] = \bb{E} \li [ \sum_{i
=2}^{\mcal{M}_{u}} Z_i - u  \ri ] \ee for all $u > 0$. Define
$\fra{M}_u = \inf \li \{ n \in \bb{N}: \sum_{i=1}^n Z_i  > u \ri \}$
for $u > 0$.  Since $Z_1, Z_2, \cd$ are i.i.d. random variables, it
must be true that $\sum_{i =2}^{\mcal{M}_{u}} Z_i$ and $\sum_{i =
1}^{\fra{M}_{u}} Z_i$ have the same distribution for all $u > 0$.
Hence, \be \la{use3698d}
 \bb{E} \li [ \sum_{i =2}^{\mcal{M}_{u}} Z_i -
u \ri ] = \bb{E} \li [ \sum_{i = 1}^{\fra{M}_{u}} Z_i - u  \ri ] \ee
for all $u > 0$.  Combining (\ref{use3698a})--(\ref{use3698d})
yields \be \la{follow88a}
 \bb{E} \li [ \li ( \sum_{i =2}^{\mscr{M}_{\bs{\lm}}}  Z_i - ( \bs{\lm} - Z_1) \ri ) \; \bb{I}_{ \{
Z_1 < \bs{\lm} \} } \ri ] = \int_{u > 0} \bb{E} \li [ \sum_{i =
1}^{\fra{M}_{u}} Z_i - u \ri ] d F_\ze (u). \ee By Lorden's
inequality \cite{Lorden}, we have \be \la{follow88b}
 \bb{E} \li [ \sum_{i = 1}^{\fra{M}_{u}} Z_i - u  \ri ] \leq \f{  \bb{E} [ (Z^+)^2 ] }{  \bb{E} [
Z ] } \ee for all $u > 0$.  Making use of (\ref{follow88a}) and
(\ref{follow88b}), we have \bel  \bb{E} \li [ \li ( \sum_{i
=2}^{\mscr{M}_{\bs{\lm}}}  Z_i - ( \bs{\lm} - Z_1) \ri ) \; \bb{I}_{
\{ Z_1 < \bs{\lm} \} } \ri ]  &   \leq  &  \int_{u > 0} \f{ \bb{E} [
(Z^+)^2 ] }{ \bb{E} [ Z ] } d F_\ze (u)
=   \f{ \bb{E} [ (Z^+)^2 ] }{ \bb{E} [ Z ] }  \int_{u > 0} d F_\ze (u) \nonumber\\
&  = & \f{ \bb{E} [ (Z^+)^2 ] }{ \bb{E} [ Z ] }  \Pr \{ \ze > 0 \}
= \f{ \bb{E} [ (Z^+)^2 ] }{ \bb{E} [ Z ] }  \Pr \{ \bs{\lm} - Z_1 > 0 \} \nonumber\\
&  =  &  \f{ \bb{E} [ (Z^+)^2 ] }{ \bb{E} [ Z ] } \Pr \{ Z <
\bs{\lm} \}. \la{lorden8899} \eel On the other hand, \be
\la{lorden8899b}
 \bb{E} [ R_{\bs{\lm}} \; \bb{I}_{ \{ Z_1 \geq \bs{\lm} \} } ] =
 \bb{E} [ (Z_1 - \bs{\lm} )^+ ] = \bb{E} [ ( Z - \bs{\lm} )^+ ]. \ee Combining (\ref{lorden8899}) and
(\ref{lorden8899b}) yields $\bb{E} [ R_{\bs{\lm}} ] = \bb{E} [
R_{\bs{\lm}} \; \bb{I}_{ \{ Z_1 < \bs{\lm} \} } ] + \bb{E} [
R_{\bs{\lm}} \; \bb{I}_{ \{ Z_1 \geq \bs{\lm} \} } ] \leq \f{ \bb{E}
[ (Z^+)^2 ]}{ \bb{E} [ Z ] } \Pr \{ Z < \bs{\lm} \} + \bb{E} [ ( Z -
\bs{\lm} )^+ ]$.  This completes the proof of the theorem.

\sect{Proof of Theorem \ref{Lordenspec} } \la{Lordenspecapp}

We need a preliminary result.

\beL \la{lem889966} Let $X_1, X_2, \cd$ be i.i.d. positive random
variables having the same distribution as $X$ such that $\bb{E} [
X^2 ] < \iy$. Define $S_n = \sum_{i=1}^n X_i$ for $n \in \bb{N}$.
Let $N_1, N_2, \cd$ be an increasing sequence of positive integers.
Define $h = \sup_{\ell \geq 0} (N_{\ell + 1} - N_\ell )$ with $N_0 =
0$. Define $\bs{N}_t = \inf \{ n \in \mscr{N}: S_N > t \}$ for $t >
0$, where $\mscr{N} = \{ N_1, N_2, \cd \}$. Define $R_t =
S_{\bs{N}_t } - t$. Then, $\bb{E} [ R_t ] \leq (h -1) \bb{E} [ X ] +
\f{ \bb{E} [ X^2 ] }{ \bb{E} [ X ] }$ for any $t > 0$.

\eeL

\bpf

Let $t > 0$.  Define $\bs{M}_t$ as the largest integer which is less
than $\bs{N}_t$ and taking value in the set $\{  N_0, N_1, N_2, \cd
\}$. Define $\bs{\mcal{N}}_t = \inf \{ n \in \bb{N}: S_n > t \}$. We
claim that $S_{h - 1 + \bs{\mcal{N}_t} } \geq S_{\bs{N}_t}$.  To
show this claim,  note that $S_k$ is increasing with respect to $k
\in \bb{N}$ as a consequence of  $X > 0$.   Since $S_{\bs{M}_t} \leq
t < S_{\bs{\mcal{N}}_t }$,  we have $\bs{M}_t \leq \bs{\mcal{N}}_t -
1$. By the definition of $h$, we have $S_{\bs{N}_t} \leq S_{h +
\bs{M}_t} \leq S_{h - 1 + \bs{\mcal{N}}_t}$.  The claim is thus
true. If follows that $\bb{E} [ S_{\bs{N}_t} ] \leq  \bb{E} [ S_{h -
1 + \bs{\mcal{N}_t} } ]$.  Since $h - 1 + \bs{\mcal{N}}_t$ is a
stopping time, by Wald's first equation, we have $\bb{E} [ S_{h - 1
+ \bs{\mcal{N}}_t} ] = ( \bb{E} [ \bs{\mcal{N}}_t ] + h - 1 ) \bb{E}
[ X ]$.  Therefore, \bee \bb{E} [ S_{\bs{N}_t} - t ] & \leq  &
\bb{E} [ S_{h - 1 + \bs{\mcal{N}_t} } - t ]
 =  \bb{E} [ S_{h - 1 + \bs{\mcal{N}}_t} - S_{\bs{\mcal{N}}_t} ] + \bb{E} [ S_{\bs{\mcal{N}}_t} - t ]\\
& = &  \bb{E} [ S_{h - 1 + \bs{\mcal{N}_t}}  ]  - \bb{E} [  S_{\bs{\mcal{N}}_t} ] + \bb{E} [ S_{\bs{\mcal{N}}_t} - t ]
 =  ( \bb{E} [ \bs{\mcal{N}}_t ] + h - 1 ) \bb{E} [ X ]
 -  \bb{E} [ \bs{\mcal{N}}_t ] \bb{E} [ X ]  + \bb{E} [ S_{\bs{\mcal{N}}_t} - t ]\\
& = & (h-1) \bb{E} [ X ]  + \bb{E} [ S_{\bs{\mcal{N}}_t} - t ]. \eee
By Lorden's inequality, $\bb{E} [ S_{\bs{\mcal{N}}_t} - t ] \leq \f{
\bb{E} [ X^2 ] }{ \bb{E} [ X ] }$.  Hence, $\bb{E} [ R_t ] = \bb{E}
[ S_{\bs{N}_t} - t ] \leq (h -1) \bb{E} [ X ]  + \f{ \bb{E} [ X^2 ]
}{ \bb{E} [ X ]  }$.  This completes the proof of the lemma.

\epf

We are now in a position to prove the theorem.  Define $\ze =
\bs{\lm} - Y$.  Let $F_\ze (.)$ denote the cumulative distribution
of $\ze$. Note that \bel  \bb{E} \li [ \li ( \sum_{i = N_1 +
1}^{\mscr{M}_{\bs{\lm}}}  Z_i - ( \bs{\lm} - Y) \ri ) \; \bb{I}_{ \{
Y < \bs{\lm} \} } \ri ]
&  =  & \bb{E} \li [ \li ( \sum_{i = N_1 + 1}^{\mscr{M}_{\bs{\lm}}}  Z_i - \ze \ri ) \; \bb{I}_{ \{ \ze > 0 \} } \ri ] \nonumber\\
& = & \int_{u > 0} \bb{E} \li [ \li ( \sum_{i = N_1 +
1}^{\mscr{M}_{\bs{\lm}}}  Z_i - \ze \ri ) \; \mid \ze = u  \ri ] d
F_\ze (u). \qqu \la{use369888a} \eel By the definition of
$\mscr{M}_{\bs{\lm}}$, we have \be \la{use369888b} \bb{E} \li [ \li
( \sum_{i = N_1 + 1}^{\mscr{M}_{\bs{\lm}}}  Z_i - \ze \ri ) \; \mid
\ze = u  \ri ] = \bb{E} \li [ \li ( \sum_{i = N_1 +
1}^{\mcal{M}_{u}} Z_i - u \ri ) \; \mid \ze = u  \ri ], \ee where
$\mcal{M}_u = \inf \li \{ n \in \mscr{N}: n \geq N_2, \; \sum_{i=
N_1 + 1}^n Z_i  > u \ri  \}$.  Since the random variables $Z_1, Z_2,
\cd$ and $\bs{\lm}$ are independent, it follows that $\ze$ and $Z_i,
\; i > N_1$ are independent.  Hence,  \be \la{use369888c}
 \bb{E} \li [ \li ( \sum_{i = N_1 + 1}^{\mcal{M}_{u}} Z_i - u \ri ) \; \mid \ze = u  \ri ] = \bb{E} \li [ \sum_{i
= N_1 + 1}^{\mcal{M}_{u}} Z_i - u  \ri ] \ee for all $u > 0$. Define
\[ \fra{M}_u = \inf \li \{ n \in \fra{N}: \sum_{i= 1}^n Z_i
> u \ri \} \qu \tx{for $u > 0$, where} \qu  \fra{N} = \{ N_\ell - N_1: \ell = 2, 3, \cd
\}.
\]
Since $Z_1, Z_2, \cd$ are i.i.d. random variables, it must be true
that $\sum_{i = N_1 + 1}^{\mcal{M}_{u}} Z_i$ and $\sum_{i =
1}^{\fra{M}_{u}} Z_i$ have the same distribution for all $u > 0$.
Hence, \be \la{use369888d}
 \bb{E} \li [  \sum_{i = N_1 + 1}^{\mcal{M}_{u}} Z_i -
u  \ri ] = \bb{E} \li [  \sum_{i = 1}^{\fra{M}_{u}} Z_i - u \ri ]
\ee for all $u > 0$.  Combining
(\ref{use369888a})--(\ref{use369888d}) yields \be \la{follow8888a}
 \bb{E} \li [ \li ( \sum_{i = N_1 + 1}^{\mscr{M}_{\bs{\lm}}}  Z_i - ( \bs{\lm} - Y) \ri ) \; \bb{I}_{ \{
Z_1 < \bs{\lm} \} } \ri ] = \int_{u > 0} \bb{E} \li [  \sum_{i =
1}^{\fra{M}_{u}} Z_i - u  \ri ] d F_\ze (u). \ee By Lemma
\ref{lem889966}, we have \be \la{follow888888b}
 \bb{E} \li [ \sum_{i = 1}^{\fra{M}_{u}} Z_i - u \ri ]
 \leq (K -1) \bb{E} [ Z ]  + \f{ \bb{E} [ Z^2 ] }{ \bb{E} [ Z ] }
  \ee for all $u > 0$.  Making use of (\ref{follow8888a}) and (\ref{follow888888b}),
  we have \bel  &   & \bb{E} \li [ \li ( \sum_{i
= N_1 + 1}^{\mscr{M}_{\bs{\lm}}}  Z_i - ( \bs{\lm} - Y) \ri ) \;
\bb{I}_{ \{ Y < \bs{\lm} \} } \ri ] \\
&   & \leq  \int_{u > 0} \li ( (K -1) \bb{E} [ Z ]  +  \f{ \bb{E} [
(Z^+)^2 ] }{ \bb{E} [ Z ] } \ri ) d F_\ze (u)  =   \li ( (K -1) \bb{E} [ Z ]
+  \f{ \bb{E} [ (Z^+)^2 ] }{ \bb{E} [ Z ] }  \ri ) \int_{u > 0} d F_\ze (u) \nonumber\\
&  & = \li ( (K -1) \bb{E} [ Z ]  +  \f{ \bb{E} [ (Z^+)^2 ] }{ \bb{E} [ Z ] }  \ri ) \Pr \{ \ze > 0 \}
 =  \li ( (K -1) \bb{E} [ Z ]  +  \f{ \bb{E} [ (Z^+)^2 ] }{ \bb{E} [ Z ] } \ri ) \Pr \{ \bs{\lm} - Y > 0 \} \nonumber\\
&   &  = \li ( (K -1) \bb{E} [ Z ]  +   \f{ \bb{E} [ (Z^+)^2 ] }{
\bb{E} [ Z ] } \ri ) \Pr \{ Y < \bs{\lm} \}. \la{lorden889988} \eel
On the other hand, \be \la{lorden889988b}
 \bb{E} [ R_{\bs{\lm}} \; \bb{I}_{ \{ Y \geq \bs{\lm} \} } ] = \bb{E} [ (Y - \bs{\lm} )^+ ]. \ee
 Combining (\ref{lorden889988}) and
(\ref{lorden889988b}) yields
\[
\bb{E} [ R_{\bs{\lm}} ] = \bb{E} [ R_{\bs{\lm}} \; \bb{I}_{ \{  Y <
\bs{\lm} \} } ] + \bb{E} [ R_{\bs{\lm}} \; \bb{I}_{ \{ Y \geq
\bs{\lm} \} } ] \leq \li ( (K -1) \bb{E} [ Z ] + \f{ \bb{E} [ Z^2 ]
}{ \bb{E} [ Z ]  } \ri ) \Pr \{  Y < \bs{\lm} \} + \bb{E} [ ( Y -
\bs{\lm} )^+ ].
\]
This completes the proof of the theorem.

\sect{Proof of Theorem  \ref{Convexlower} } \la{Convexlowerapp}

We shall first show $\bb{E} [ \bs{N} ] \geq \mscr{A} (\mu)$ under
the assumption that $\mscr{A} (\mu) < \iy$.  If $\bb{E} [ \bs{N} ] =
\iy$, then $\bb{E} [ \bs{N} ] \geq \mscr{A} (\mu)$ trivially holds.
If $\bb{E} [ \bs{N} ] < \iy$, then $\Pr \{ \bs{N} < \iy  \} = 1$ and
it follows that $S_{\bs{N}}$ is well-defined and $\Pr \{ (\bs{N},
S_{\bs{N}}) \in \mscr{R}^c \} = 1$.  Since $\bb{E} [ \bs{N} ] <
\iy$, it follows from Wald's equation that $\bb{E} [  S_{\bs{N}}] =
\bb{E} [ \bs{N} ] \mu$.   According to Theorem \ref{fundamental}, we
have $( \bb{E} [ \bs{N} ], \;  \bb{E} [ S_{\bs{N}}] ) \in
\mscr{R}^c$.  Using Wald's equation, we have $( \bb{E} [ \bs{N} ],
\; \bb{E} [ \bs{N} ] \mu ) \in \mscr{R}^c$. It follows from the
definition of IDET that $\bb{E} [ \bs{N} ] \geq \mscr{A} (\mu)$.

It remains to show $\bb{E} [ \bs{N} ] = \iy$ under the assumption
that $\mscr{A} (\mu) = \iy$.  We use a contradiction method. Suppose
that $\bb{E} [ \bs{N} ] < \iy$, then $\Pr \{ \bs{N} < \iy  \} = 1$
and it follows that $\Pr \{ (\bs{N}, S_{\bs{N}}) \in \mscr{R}^c \} =
1$. Since $\bb{E} [ \bs{N} ] < \iy$, it follows from Wald's equation
that $\bb{E} [ S_{\bs{N}}] = \bb{E} [ \bs{N} ] \mu$.  According to
Theorem \ref{fundamental}, we have $( \bb{E} [ \bs{N} ], \;  \bb{E}
[  S_{\bs{N}}] ) \in \mscr{R}^c$.  Applying Wald's equation, we have
$( \bb{E} [ \bs{N} ], \; \bb{E} [ \bs{N} ] \mu  ) \in \mscr{R}^c$.
It follows from the definition of IDET that  $\mscr{A} (\mu) \leq
\bb{E} [ \bs{N} ] < \iy$.  This is a contradiction.  Therefore, it
must be true that $\bb{E} [ \bs{N} ] = \iy$ if $\mscr{A} (\mu) =
\iy$.  The proof of the theorem is thus completed.

\sect{Proof of Theorem \ref{BoundGen} } \la{BoundGenapp}

We need a preliminary result.

\beL \la{Lem33669988}

 There exist a row matrix $A$ of size $1 \times d$,  a real number $B$,
 and a positive number $C$ such that $A \mu + B > 0$ and
 that $\mscr{R} \subseteq \{  (t, s): t \geq 0, \; s \in \bb{R}^d,    As + B t < C \}$.

\eeL

\bpf

For simplicity of notations, define $m = \mscr{B} (\mu)$.  Let $\vep
> 0$.  Recall that $\mscr{B} (\mu) = \sup \{ t \geq 0: (t, \mu t) \in \mscr{R} \}$.
By the assumption that $\mscr{B} (\mu) < \iy$, it follows that $(m,
(m + \vep) \mu )$ must be an interior point of the complementary set
$\mscr{R}^c$.  Since $\mscr{R}$ is convex, it follows from the {\it
hyperplane separation theorem} \cite[Theorem 4.11, page
170]{Bertsimas} that there exists a hyperplane which strictly
separate the point $(m+ \vep, (m + \vep) \mu )$ and the convex set
$\mscr{R}$.  This implies that there exist a row matrix $A$ of size
$1 \times d$,   real numbers $B$ and $C$ such that $A s + B t < C$
for $(t, s) \in \mscr{R}$ and that $A s + B t > C$ for $(t, s) = (m
+ \vep, (m + \vep) \mu )$.  Since $\mscr{R}$ contains $(0,
\bs{0}_d)$, it must be true that $C > 0$. Since $\mscr{R}$ contains
$ (m, m \mu)$, it follows that $m (A \mu + B) < C$.
 Since $\mscr{R}$ does not contain $ (m + \vep, (m + \vep) \mu)$, it follows
that $(m + \vep) (A \mu + B) > C$.  Hence, $m (A \mu + B) < C < (m +
\vep) (A \mu + B)$. This implies that $\vep (A \mu + B) > 0$.  Since
$\vep > 0$, it follows that $A \mu + B > 0$. This completes the
proof of the lemma.

\epf

\bsk

We are now in a position to prove the theorem.  Let $A, \; B$ and
$C$ be defined as in Lemma \ref{Lem33669988}.  Define
\[
Y = A X + B, \qqu Y_i = A X_i + B, \qqu i \in \bb{N}.
\]
It follows from Lemma \ref{Lem33669988} that $\bb{E} [ Y ] = \bb{E}
[ A X + B ] = A \mu + B > 0$.  By the dominant convergence theorem,
$\lim_{n \to \iy} \bb{E} [ Y \bb{I}_{ \{ Y \leq n \} } ] = \bb{E} [
Y ] > 0$.  Therefore, there exists a positive number $M$ such that
$\bb{E} [ Y \bb{I}_{ \{  Y \leq M   \} } ] > 0$.  Define
\[
Z = Y \bb{I}_{ \{  Y \leq M   \} }, \qqu Z_i = Y_i \bb{I}_{ \{  Y_i
\leq M   \} }, \qqu i \in \bb{N}
\]
and $\mscr{S}_n = \sum_{i=1}^n Z_i$ for $n \in \bb{N}$.  Clearly,
$Z_1, Z_2, \cd$ are i.i.d. random variables having the same
distribution as $Z$ with $\Pr \{ Z \leq M  \} = 1$ and $\bb{E} [ Z ]
> 0$.  Define
\[
\bs{T} = \{  n \in \mscr{N}:    A \mscr{S}_n + B n > C    \}.
\]
Note that $A \mscr{S}_n + B n = \sum_{i = 1}^n Z_i \leq \sum_{i =
1}^n Y_i = A S_n + B n$ for $n \in \bb{N}$. By Lemma
\ref{Lem33669988}, we have  $\mscr{R} \subseteq \{  (t, s): t \geq
0, \; s \in \bb{R}^d,    As + B t < C \}$.   This implies that
\[
\{ A \mscr{S}_n + B n > C \}  \subseteq \{ A S_n + B n > C \}
\subseteq \{ (n, S_n) \notin \mscr{R} \}
\]
for all $n \in \bb{N}$.  Therefore, $\bs{N} \leq \bs{T}$.  Hence, to
show $\bb{E} [ \bs{N} ] < \iy$,  it suffices to show $\bb{E} [
\bs{T} ] < \iy$. Define
\[ \bs{T}_n = \min \{ \bs{T}, \;  n \},  \qqu n \in \mscr{N}.
\]
Then, $\bs{T}_n$ is a stopping time and $\bb{E} [ \bs{T}_n ] < \iy$.
By Wald's equation, $\bb{E} [ \mscr{S}_{ \bs{T}_n } ] = \bb{E} [ Z ]
\bb{E} [ \bs{T}_n ]$ for $n \in \mscr{N}$.  Let $\eta$ be a number
satisfying \be \la{lmchoice} 0 < \eta < \min \li \{ \f{1}{2}, \;
\f{\bb{E} [ Z ] } {  M } \ri \}. \ee As a consequence of the
assumption that $\lim_{\ell \to \iy} \f{N_{\ell + 1}}{N_\ell} = 1$,
there exists a number $\ell^*$ such that $N_\ell - N_{\ell - 1}  <
\eta N_\ell$ for all $\ell \geq \ell^*$.  We claim that
\[
\mscr{S}_{ \bs{T}_n } \leq  C + M \eta (N_{\ell^*} + \bs{T}_n) \qqu
\tx{for all $n \in \mscr{N}$.}
\]
To prove the claim, we proceed as follows.

Let $\om \in \Om$.   If $\bs{T} (\om) = \iy$, then the claim holds
trivially.  Hence, it suffices to consider the scenarios that
$\bs{T} (\om) = N_\ell$ for some $\ell \in \bb{N}$. There are four
cases:

Case (i):  $\ell \geq \ell^*, \; n \geq N_\ell$;

Case (ii):  $\ell \geq \ell^*, \; n < N_\ell$;

Case (iii):  $\ell < \ell^*, \; n \geq N_\ell$;

Case (iv):  $\ell < \ell^*, \; n < N_\ell$.

\bsk

In Case (i), we have $\bs{T}_n (\om)  = \min \{ N_\ell, \; n \} =
N_\ell = \bs{T} (\om)$ and $\mscr{S}_{ \bs{T}_n } (\om) = \mscr{S}_{
N_\ell } (\om)
> C \geq \mscr{S}_{ N_{\ell - 1} } (\om)$.  Since $N_\ell - N_{\ell - 1} \leq \eta N_\ell$ for
$\ell \geq \ell^*$, it follows that $N_\ell - N_{\ell - 1} < \eta
(N_{\ell^*} + N_\ell)$ for $\ell \geq \ell^*$. So,
\[
\mscr{S}_{ N_\ell } (\om) = \mscr{S}_{ N_{\ell - 1} } (\om) +
\mscr{S}_{ N_\ell } (\om) - \mscr{S}_{ N_{\ell - 1} } (\om) \leq C +
(N_\ell - N_{\ell - 1}) M < C + \eta (N_{\ell^*} + N_\ell) M.
\]

In Case (ii), we have $\bs{T}_n (\om)  = \min \{ N_\ell, \; n \} = n
\leq N_{\ell -1}$.  Hence, $\mscr{S}_{ \bs{T}_n } (\om) = \mscr{S}_{
n } (\om) \leq C$ because $n \leq N_{\ell -1}$.

In Case (iii), we have $\bs{T}_n  (\om) = \min \{ N_\ell, \; n \} =
N_\ell = \bs{T} (\om)$ and $\mscr{S}_{ \bs{T}_n } (\om) = \mscr{S}_{
N_\ell } (\om) > C \geq \mscr{S}_{ N_{\ell - 1} } (\om)$.  Since
$\eta \in (0, \f{1}{2})$ and the sequence $\{ N_\ell \}$ is
increasing with respect to $\ell$, it follows that $N_\ell \leq
\f{\eta}{1 - \eta} N_{\ell^*}$ for $\ell < \ell^*$. Hence, $N_\ell
\leq \eta (N_{\ell^*} + N_\ell)$ for $\ell < \ell^*$. Of course,
$N_\ell - N_{\ell - 1} \leq \eta (N_{\ell^*} + N_\ell)$ for $\ell <
\ell^*$.  So,
\[
\mscr{S}_{ \bs{T}_n } (\om) = \mscr{S}_{ N_\ell } (\om) = \mscr{S}_{
N_{\ell - 1} } (\om) + \mscr{S}_{ N_\ell } (\om) - \mscr{S}_{
N_{\ell - 1} } (\om) \leq C + (N_\ell - N_{\ell - 1}) M \leq C +
\eta (N_{\ell^*} + N_\ell) M.
\]

In Case (iv), we have $\bs{T}_n  (\om) = \min \{ N_\ell, \; n \} = n
\leq N_{\ell -1}$.  Hence, $\mscr{S}_{ \bs{T}_n } (\om) = \mscr{S}_{
n } (\om) \leq C$ because $n \leq N_{\ell -1}$.

\bsk

Therefore, we have $\mscr{S}_{ \bs{T}_n } (\om) \leq  C + M \eta [
N_{\ell^*} + \bs{T}_n (\om) ]$ for all cases. This proves the claim.

Since $\mscr{S}_{ \bs{T}_n } \leq  C + M \eta (N_{\ell^*} +
\bs{T}_n)$ for all $n \in \mscr{N}$,  taking expectations on both
sides of this inequality and applying Wald's equation  yields
\[
( \bb{E} [ Z ] - \eta M) \bb{E} [ \bs{T}_n ] \leq C + \eta M
N_{\ell^*}.
\]
From (\ref{lmchoice}), we have $\bb{E} [ Z ] > \eta M$.  It follows
that
\[
\bb{E} [ \bs{T}_n ] \leq \f{  C + \eta M  N_{\ell^*} }{  \bb{E} [ Z
] - \eta M} < \iy.
\]
Note that $\{ \bs{T}_n, \; n \in \mscr{N} \}$ is a sequence of
positive random variables convergent to $\bs{T}$ as $n \to \iy$.  By
Fatou's lemma,
\[
\bb{E} [ \bs{T} ] \leq  \liminf_{n \to \iy} \bb{E} [ \bs{T}_n ] \leq
\f{  C + \eta M  N_{\ell^*} }{  \bb{E} [ Z ]  - \eta M} < \iy.
\]
It follows that $\bb{E} [ \bs{N} ] \leq \bb{E} [ \bs{T} ] < \iy$.
This completes the proof of the theorem.

\sect{Proof of Theorem \ref{GenConvex} } \la{GenConvexapp}

Since conditions (I)--(IV) are fulfilled, it follows from Theorem
\ref{BoundGen} that $\bb{E} [ \bs{M} ] < \bb{E} [ \bs{N} ] < \iy$,
which implies $\Pr \{ \bs{M} < \iy \} = 1 $ and $\Pr \{ \bs{N} < \iy
\} = 1$. Hence, $\bs{N}$ and $\bs{M}$ are well-defined random
variables.  Define
\[
\vDe = S_{\bs{N}} - S_{\bs{M}} - (\bs{N} - \bs{M}) \mu. \]
Our proof of the theorem relies on some properties of $\vDe$ as stated by the
following lemma.

\beL \la{ineq9996688}
\bel &  &  \bb{E} [ \vDe^+ ] \leq \f{1}{2} \bb{E} [ \bs{N} - N_0 ] \; \xi, \la{expvdep}\\
&  &  \bb{E} [ \vDe^- ] \leq \f{1}{2} \bb{E} [ \bs{N} - N_0 ] \; \xi
\la{expvden},   \eel where $\xi = \bb{E} [ |X - \mu|]$. \eeL

\bpf
 Define $\vDe_\ell = S_{N_\ell} - S_{N_{\ell - 1}} - (N_\ell - N_{\ell - 1})
\mu$
 for $\ell \in \bb{N}$.  Let $\bs{\tau}$ denote the stopping index such that $N_{\bs{\tau}} = \bs{N}$. Note that
 \bee
\bb{E} [ \vDe^+ ] & = & \sum_{\ell = 1}^\iy \bb{E} [ \vDe^+ \; \bb{I}_{ \{ \bs{\tau} = \ell \} } ]  =
 \sum_{\ell = 1}^\iy \bb{E} [ \vDe_\ell^+ \; \bb{I}_{ \{ \bs{\tau} = \ell \} } ]
 =  \bb{E} [ \vDe_1^+ \;  \bb{I}_{ \{ \bs{\tau} = 1 \} } ]
 + \sum_{\ell = 2}^\iy \bb{E} [ \vDe_\ell^+ \; \bb{I}_{ \{ \bs{\tau} = \ell \} } ]\\
& \leq & \bb{E} [ \vDe_1^+ \; \bb{I}_{ \{ \bs{\tau} = 1 \} } ] +
\sum_{\ell = 2}^\iy \bb{E} [ \vDe_\ell^+ \; \bb{I}_{ \{ \bs{\tau} >
\ell - 1 \} } ] \leq  \bb{E} [ \vDe_1^+] + \sum_{\ell = 2}^\iy
\bb{E} [ \vDe_\ell^+ \; \bb{I}_{ \{ \bs{\tau} > \ell - 1 \} } ].
 \eee
Observing that $\vDe_\ell$ depends only on
$\{ X_n: N_{\ell - 1} + 1 \leq n \leq N_\ell \}$ and that the event $\{ \bs{\tau} > \ell - 1 \}$
depends only on $\{ X_n: 1 \leq n \leq N_{\ell - 1} \}$, we have that
\[
\bb{E} [ \vDe_\ell^+ \; \bb{I}_{ \{ \bs{\tau} > \ell - 1 \} } ] =
\bb{E} [ \vDe_\ell^+ ] \; \bb{E} [ \bb{I}_{ \{ \bs{\tau} > \ell - 1
\} } ]  = \bb{E} [ \vDe_\ell^+ ] \Pr \{  \bs{\tau} > \ell - 1 \}
\]
for $\ell > 1$.  It follows that
\[
\bb{E} [ \vDe^+ ] \leq \bb{E} [ \vDe_1^+]
+ \sum_{\ell = 2}^\iy \bb{E} [ \vDe_\ell^+ \; \bb{I}_{ \{ \bs{\tau} > \ell - 1 \} } ]  =  \bb{E} [
\vDe_1^+] + \sum_{\ell = 2}^\iy \bb{E} [ \vDe_\ell^+ ] \Pr \{ \bs{\tau} > \ell - 1 \}.
\]
Since $X_1, X_2, \cd$ are random vectors having the same
distribution as $X$, we have
 \[
\bb{E} [ \vDe_\ell^+ ] \leq (N_\ell - N_{\ell - 1}) \bb{E} [  (X -
\mu)^+] = \f{1}{2} (N_\ell - N_{\ell - 1}) \xi, \qqu \ell \in
\bb{N}.
 \]
Hence,
 \bee
\bb{E} [ \vDe^+ ]  \leq  \f{1}{2}  \li [ N_1 - N_0 + \sum_{\ell =
1}^\iy (N_{\ell+ 1} - N_{\ell})  \Pr \{ \bs{\tau} > \ell \} \ri ]
\xi = \f{1}{2} \bb{E} [ \bs{N} - N_0 ] \; \xi.
 \eee
 This proves (\ref{expvdep}).
 By similar arguments we can show the inequalities (\ref{expvden})
 regarding $\bb{E} [ \vDe^- ]$.

\epf

\beL

\la{lem1088}   Let $r = N_0 - K$ and $\xi = \bb{E} [  | X - \mu |]$.
Then, $\li | \bb{E} [ S_{\bs{M}} ]  - \bb{E} [ \bs{M} ] \mu \ri |
\leq \f{\lm}{2} \bb{E} [ \bs{M} ] \xi  - \f{r}{2} \xi$.  \eeL

\bpf

By the assumption that $\bb{E} [ | X | ]$ is bounded, we have that
both $\bb{E} [ X^+ ]$ and $\bb{E} [ X^- ]$ are bounded. By Theorem
\ref{BoundGen}, $\bs{N}$ is a stopping time such that  $\bb{E} [
\bs{N} ] < \iy$.  Hence, it follows from Wald's first equation that
\[ \bb{E} [ (S_{\bs{N}})^+ ] \leq \sum_{i=1}^{\bs{N}} (X_i)^+ =
\bb{E} [  \bs{N} ] \bb{E} [ X^+ ] < \iy, \qqu \bb{E} [
(S_{\bs{N}})^- ] \leq \sum_{i=1}^{\bs{N}} (X_i)^- = \bb{E} [  \bs{N}
] \bb{E} [ X^- ] < \iy. \] Thus, \be \la{finite1}
 \bb{E} [ | S_{\bs{N}} | ] \leq \max \{ \bb{E} [ (S_{\bs{N}})^+ ], \bb{E} [ (S_{\bs{N}})^- ]  \} < \iy. \ee
 By the definition of $\vDe $, we have
\be \la{basic8896}
 S_{\bs{M}} = S_{\bs{N}} - \vDe - (\bs{N} - \bs{M}) \mu.
\ee Hence, \be \la{finite2}
 \bb{E} [ | S_{\bs{M}} | ] \leq \bb{E} [ | S_{\bs{N}} | ] + \bb{E} [ |\vDe| ]  + \bb{E} [ \bs{N} - \bs{M} ] |\mu|.
\ee From  Lemma \ref{ineq9996688}, we have \be \la{finite3}
 \bb{E} [ | \vDe | ] = \bb{E} [ \vDe^+ ] + \bb{E} [ \vDe^- ] \leq \bb{E} [ \bs{N} - N_0 ] \; \xi < \iy. \ee Since
$\bb{E} [ \bs{M} ] < \bb{E} [ \bs{N} ] < \iy$, we have \be
\la{finite4}
 \bb{E} [ \bs{N} - \bs{M} ] < \iy. \ee
 Combining (\ref{finite1})--(\ref{finite4}) leads to
 \be
 \la{boundSM88}
 \bb{E} [ | S_{\bs{M}} | ] < \iy.
 \ee
This establishes the existence of $\bb{E} [ S_{\bs{M}} ]$.  Taking
expectations on both sides of (\ref{basic8896}) yields
\bel \bb{E} [ S_{\bs{M}} ] & = & \bb{E} [ S_{\bs{N}} - \vDe - (\bs{N} - \bs{M}) \mu ] \nonumber\\
& = & \bb{E} [ S_{\bs{N}} ] - \bb{E} [  \vDe ] - \bb{E} [  \bs{N} ] \mu + \bb{E} [  \bs{M} ] \mu \nonumber\\
& = & \bb{E} [  \bs{N} ] \mu - \bb{E} [  \vDe ] - \bb{E} [  \bs{N} ] \mu + \bb{E} [  \bs{M} ] \mu \la{usewald} \\
& = & \bb{E} [  \bs{M} ] \mu  - \bb{E} [  \vDe ] \nonumber\\
& = &  \bb{E} [  \bs{M} ] \mu  - \bb{E} [  \vDe^+ ] +  \bb{E} [  \vDe^- ], \la{touse88}  \eel
where we have used Wald's equation $\bb{E} [
S_{\bs{N}} ] = \bb{E} [  \bs{N} ] \mu$ in (\ref{usewald}).   As a consequence of (\ref{touse88}), we have \be \la{ineqvip998}
 \bb{E} [  \bs{M} ]
\mu - \bb{E} [  \vDe^+ ] \leq \bb{E} [ S_{\bs{M}} ]  \leq \bb{E} [
\bs{M} ] \mu  + \bb{E} [  \vDe^- ]. \ee In view of $\bs{N} \leq  \lm
\bs{M} + K$, we have \be \la{vipineq} \bb{E} [ \bs{N} - N_0 ] \leq
\lm \bb{E} [ \bs{M} ] - r. \ee  Making use of (\ref{expvden}),
(\ref{vipineq}) and the second inequality of (\ref{ineqvip998}), we
have  \bee  \bb{E} [ S_{\bs{M}} ]  \leq  \bb{E} [ \bs{M} ] \mu +
\bb{E} [ \vDe^- ] \leq  \bb{E} [ \bs{M} ] \mu + \f{1}{2} \bb{E} [
\bs{N} - N_0 ] \; \xi  \leq  \bb{E} [ \bs{M} ] \mu + \f{\lm}{2}
\bb{E} [ \bs{M} ] \xi  - \f{r}{2}  \xi.  \eee Making use of
(\ref{expvdep}),  (\ref{vipineq}) and first inequality of
(\ref{ineqvip998}), we have \bee \bb{E} [ S_{\bs{M}} ]  \geq \bb{E}
[ \bs{M} ] \mu - \bb{E} [ \vDe^+ ] \geq  \bb{E} [ \bs{M} ] \mu -
\f{1}{2} \bb{E} [ \bs{N} - N_0 ] \; \xi \geq  \bb{E} [ \bs{M} ] \mu
- \f{\lm}{2} \bb{E} [ \bs{M} ] \xi  + \f{r}{2}  \xi.  \eee This
completes the proof of the lemma.  \epf

We are now in a position to prove the theorem.  By the definition of
$\bs{M}$, we have $\Pr \{ (\bs{M}, S_{\bs{M}} ) \in \mscr{R} \} =
1$.  Since $\bb{E} [ \bs{M} ] < \iy$ and $\bb{E} [ S_{\bs{M}} ]$
exists, it follows from Theorem \ref{fundamental} that $( \bb{E} [
\bs{M} ], \bb{E} [ S_{\bs{M}} ] ) \in \mscr{R}$.  The conclusion of
the theorem immediately follows from this fact and Lemma
\ref{lem1088}.

\sect{Proof of Theorem \ref{Convexsecond}} \la{Convexsecondapp}

We need some preliminary results.

\beL \la{convexupper}

 Let $Z$ be a scalar random variable with mean $\se$
 such that $\Pr \{ a \leq Z \leq b \} = 1$, where $a < b$ are real numbers.
  If $g(x)$ is a convex function of $x \in [a,
b]$, then \be \la{chain}
 \bb{E}[ g (Z) ] \leq \f{1}{b - a} [ (b - \se) g(a) +
( \se - a) g(b) ]. \ee In particular,  \bel &  & \bb{E} [ ( Z - \se
)^+ ] = \bb{E} [ ( Z - \se )^- ] = \f{1}{2} \bb{E} [ | Z - \se | ]
\leq  \f{ (\se - a) (b - \se) }{b - a}. \la{app89b} \la{app89c} \eel

\eeL

\bpf

To show (\ref{chain}), note that, as a consequence of the convexity of the function $g$,
\[
g(x) \leq \f{ g(b) - g(a) }{b - a}  (x - a) + g(a), \qqu x \in [a, b].
\]
By the assumption that $\Pr \{ a \leq Z \leq b \} = 1$, we have \[
g(Z) \leq \f{ g(b) - g(a) }{b - a}  (Z - a) + g(a)
\]
almost surely.  Taking expectation on both sides of the above inequality yields
\[
\bb{E}[ g (Z) ] \leq \f{ g(b) - g(a) }{b - a} \bb{E} [ Z - a ] +
g(a) = \f{1}{b - a} [ (b - \se) g(a) + (\se - a) g(b) ].
\]
This establishes (\ref{chain}).  Applying (\ref{chain}) to convex
function $g (x) = | x - \se |$ yields (\ref{app89b}).

 \epf

 \beL  \la{multibound}
 Assume that $\Pr \{ \bs{a} \leq X \leq \bs{b} \} = 1$,
 where $\bs{a}, \; \bs{b} \in \bb{R}^d$, and that conditions (I)--(III), V and (VI) are fulfilled.
Define $\bs{v} = \f{ (\mu - \bs{a})  (\bs{b} - \mu)}{ \bs{b} -
\bs{a} }$.  Then, \bel &   & | \bb{E} [ S_{\bs{M}} ] - \bb{E} [ \bs{M} ] \mu |
 \leq \lm \bb{E} [ \bs{M}  ] \bs{v} - r \bs{v},  \la{last8899a}\\
&  &  \li \{  (\lm - 1) \bb{E} [ \bs{M} ] + K \ri \}  ( \mu - \bs{b}
)  \leq \bb{E} [ S_{\bs{M}} ] - \bb{E} [ \bs{M} ] \mu \leq  \li \{
(\lm - 1) \bb{E} [ \bs{M} ] + K \ri \}  ( \mu - \bs{a} ).
\la{last8899} \eel

\eeL

\bpf  Since $\Pr \{ \bs{a} \leq X \leq \bs{b} \} = 1$ and conditions
(I)--(III) are fulfilled, it follows from Theorem \ref{BoundGen}
that $\bb{E} [ \bs{M} ] <  \bb{E} [ \bs{N} ] < \iy$. Since $\Pr \{
\bs{a} \leq X \leq \bs{b} \} = 1$, it follows from Lemma
\ref{convexupper} that \be \bb{E} [ ( X - \mu )^+ ] = \bb{E} [ ( X -
\mu )^- ]  = \f{1}{2} \xi \leq  \bs{v}. \la{app89b88} \ee Making use
of (\ref{app89b88}) and Lemma \ref{ineq9996688}, we have \bel & &
\bb{E} [ \vDe^+ ] \leq
\bb{E} [ \bs{N} - N_0 ] \; \bb{E} [  (X - \mu)^+] \leq \bb{E} [ \bs{N} - N_0 ] \bs{v}, \la{expvdep88}\\
&  &  \bb{E} [ \vDe^- ] \leq \bb{E} [ \bs{N} - N_0 ] \; \bb{E} [  (X
- \mu)^-] \leq \bb{E} [ \bs{N} - N_0 ] \bs{v}.  \la{expvden88} \eel
Making use of (\ref{ineqvip998}), (\ref{expvdep88}),
(\ref{expvden88}) and the fact that $\bs{N} \leq \lm \bs{M} + K$, we
have \bee  \bb{E} [ S_{\bs{M}} ]  \leq  \bb{E} [  \bs{M} ] \mu  +
\bb{E} [  \vDe^- ]  \leq  \bb{E} [ \bs{M} ] \mu + \bb{E} [ \bs{N} -
N_0 ] \bs{v} \leq \bb{E} [ \bs{M} ] \mu + \lm \bb{E} [ \bs{M}  ]
\bs{v} - r \bs{v}  \eee and \bee  \bb{E} [ S_{\bs{M}} ]  \geq \bb{E}
[ \bs{M} ] \mu - \bb{E} [  \vDe^+ ] \geq  \bb{E} [ \bs{M} ] \mu -
\bb{E} [ \bs{N} - N_0 ] \bs{v} \geq  \bb{E} [ \bs{M} ] \mu - \lm
\bb{E} [ \bs{M}  ] \bs{v} + r \bs{v}.   \eee This proves
(\ref{last8899a}). It remains to show (\ref{last8899}).  Recall that
\[
\vDe = S_{\bs{N}} - S_{\bs{M}} - (\bs{N} - \bs{M} ) \mu = \sum_{i= \bs{M} + 1}^{\bs{N}} (X_i - \mu).
\]
Since $\Pr \{ \bs{a} \leq X \leq \bs{b} \} = 1$, it follows that $(X
- \mu)^+ \leq \bs{b} - \mu$ and $(X - \mu)^- \leq \mu - \bs{a}$
almost surely.  Hence,  \be \la{mak89a}
 \bb{E} [ \vDe^+ ] \leq \bb{E} \li [  \sum_{i= \bs{M} + 1}^{\bs{N}} (X_i - \mu)^+  \ri ] \leq \bb{E} [ \bs{N} - \bs{M} ]
( \bs{b} - \mu), \ee \be \la{mak89b} \bb{E} [ \vDe^- ]
\leq \bb{E} \li [  \sum_{i= \bs{M} + 1}^{\bs{N}} (X_i - \mu)^-  \ri ] \leq \bb{E} [
\bs{N} - \bs{M} ] (\mu - \bs{a}). \ee
Making use of (\ref{ineqvip998}), (\ref{mak89a}), (\ref{mak89b}) and the fact that $\bs{N} \leq  \lm
\bs{M} + K$, we have \bee  \bb{E} [ S_{\bs{M}} ]  \leq  \bb{E} [  \bs{M} ] \mu  + \bb{E} [  \vDe^- ]
 \leq  \bb{E} [ \bs{M} ] \mu + \bb{E} [ \bs{N} - \bs{M} ] ( \mu - \bs{a} )
 \leq  \bb{E} [ \bs{M} ] \mu + \li \{  (\lm - 1) \bb{E} [ \bs{M} ]
+ K \ri \}  ( \mu - \bs{a} ),  \eee and
\bee  \bb{E} [ S_{\bs{M}} ] \geq  \bb{E} [ \bs{M} ] \mu - \bb{E} [  \vDe^+ ]
 \geq  \bb{E} [ \bs{M} ] \mu - \bb{E} [ \bs{N} - \bs{M} ] ( \bs{b} - \mu)
 \geq  \bb{E} [ \bs{M} ] \mu - \li \{  (\lm - 1) \bb{E} [ \bs{M} ]
+ K \ri \}  ( \bs{b} - \mu ).   \eee This proves (\ref{last8899}).
The proof of the lemma is thus completed.

\epf

We are now in a position to prove the theorem.  By the definition of
$\bs{M}$, we have $\Pr \{ (\bs{M}, S_{\bs{M}} ) \in \mscr{R} \} =
1$. Since $\bb{E} [ \bs{M} ] < \iy$ and $\bb{E} [ S_{\bs{M}} ]$
exists as asserted by (\ref{boundSM88}), it follows from Theorem
\ref{fundamental} that $( \bb{E} [ \bs{M} ], \bb{E} [ S_{\bs{M}} ] )
\in \mscr{R}$. The conclusion of the theorem immediately follows
from this fact and Lemma \ref{multibound}.

\sect{Proof of Theorem  \ref{ConvexThird}} \la{ConvexThirdapp}

We need a preliminary result.

\beL

\la{lem88669933}

$\bb{E} [ \bs{M} ] \leq g \li (  \f{ \bb{E} [ S_{\bs{M}} ] }{ \bb{E}
[ \bs{M} ] } \ri )$.

\eeL

\bpf

Note that $\bs{N} = \{ n \in \mscr{N}: (n, S_n) \notin \mscr{R} \}$,
where $\mscr{R} = \{ (t, s): s\in \bb{R}^d, \f{s}{t} \in D,  \; 0 <
t \leq g (\f{s}{t}) \}$.  By Lemma \ref{BoundaryCon}, $t g
(\f{s}{t})$ is a concave function for $t > 0, \; s \in \bb{R}^d$
such that $\f{s}{t} \in D$. Hence, $\mscr{R} = \{ (t, s): s\in
\bb{R}^d, \f{s}{t} \in D,  \; 0 < t \leq g (\f{s}{t}) \} = \{ (t,
s): s\in \bb{R}^d, \; t > 0, \; \f{s}{t} \in D, \; 0 < t^2 \leq t g
(\f{s}{t}) \}$ is a convex set.  Note that
\[
\mscr{B} (v) = \sup \{ t > 0: (t, v t) \in \mscr{R} \} = g (v)
\]
for $v \in D$.  From Theorem \ref{BoundGen}, we have  $\bb{E} [
\bs{M} ] < \bb{E} [ \bs{N} ] < \iy$, it follows that $S_{\bs{M}}$ is
well-defined.  By the definition of the stopping rule, we have $\Pr
\{ ( \bs{M}, S_{\bs{M}} ) \in  \mscr{R} \} = 1$.  As in the proof of
Lemma \ref{lem1088}, we have that $\bb{E} [ S_{\bs{M}} ]$ exists. By
Theorem \ref{fundamental}, we have $( \bb{E} [ \bs{M} ] , \bb{E} [
S_{\bs{M}} ] ) \in  \mscr{R}$,  i.e.,
\[
\li ( \bb{E} [ \bs{M} ] , \f{ \bb{E} [ S_{\bs{M}} ] }{  \bb{E} [
\bs{M} ] } \bb{E} [ \bs{M} ] \ri ) \in \mscr{R}.
\]
Hence, $\bb{E} [ \bs{M} ] \leq \mscr{B} \li (  \f{ \bb{E} [
S_{\bs{M}} ] }{ \bb{E} [ \bs{M} ] } \ri ) = g \li (  \f{ \bb{E} [
S_{\bs{M}} ] }{ \bb{E} [ \bs{M} ] } \ri )$.  This completes the
proof of the lemma.

\epf

We are now in a position to prove the theorem.  By Lemma
\ref{lem1088}, $\li | \bb{E} [ S_{\bs{M}} ]  - \bb{E} [ \bs{M} ] \mu
\ri | \leq \f{\lm}{2} \bb{E} [ \bs{M} ] \xi  - \f{r}{2} \xi$, where
$r = N_0 - K$. As a consequence of $N_{\ell + 1} - N_\ell \leq K
\leq N_0$ for $\ell \in \bb{Z}^+$, we have that $r \leq 0$ and it
follows that $\li | \bb{E} [ S_{\bs{M}} ]  - \bb{E} [ \bs{M} ] \mu
\ri | \leq \f{\lm}{2} \bb{E} [ \bs{M} ] \xi$.  Hence, $\li | \f{
\bb{E} [ S_{\bs{M}} ] }{ \bb{E} [ \bs{M} ] } -  \mu  \ri | \leq
\f{1}{2}
 \xi$. Finally, the conclusion of the theorem follows from this inequality and Lemma \ref{lem88669933}.

\sect{Proof of Corollary \ref{Convexfour} } \la{Convexfourapp}

Define
\[
\bs{\mcal{N}} (\vep) = \inf \{ n \in \bb{N}: n >  1 + \vep +
g(\ovl{X}_n) \},
\]
where $\vep > 0$.  Define $K = 1, \; N_0 = 1$ and $N_\ell = \ell +
1$ for $\ell \in \bb{N}$. Then, $N_{\ell + 1} - N_\ell \leq K \leq
N_0$ for $\ell \in \bb{Z}^+$. Since $g$ is non-negative, it must be
true that $N_0 \leq 1 + \vep + g (\ovl{X}_{N_0})$ is a sure event.
Therefore, we can apply Theorem  \ref{ConvexThird} to stopping time
$\bs{\mcal{N}} (\vep)$ to conclude that
\[
\bb{E} [ \bs{\mcal{N}} (\vep) ]  \leq
1 + \max_{\se \in \mscr{D}} [ 1 + \vep + g(\se) ] \leq 2 + \vep + \max_{\se \in \mscr{D}} g (\se).
\]
Since $\bs{N} \leq \bs{\mcal{N}} (\vep)$, we have
\[
\bb{E} [ \bs{N} ] \leq \bb{E} [ \bs{\mcal{N}} (\vep) ] \leq 2 + \vep + \max_{\se \in \mscr{D}} g (\se).
\]
Since the above inequalities hold for arbitrarily small $\vep > 0$,
it must be true that $\bb{E} [ \bs{N} ] \leq 2 + \max_{\se \in
\mscr{D}} g (\se)$.  This completes the proof of the corollary.

\sect{Proof of Corollary  \ref{try88} } \la{try88app}

As a consequence of the assumption that $\Pr \{ \bs{a} \leq X \leq
\bs{b} \} = 1$, it must be true that each element of $\bb{E} [ | X |
]$ is finite.  It follows from Theorem \ref{ConvexThird} that
$\bb{E} [ \bs{N} ] \leq K + \max_{\se \in \mscr{D} } g (\se)$,
where $\mscr{D} = \{ \se \in D: | \se - \mu | \leq \f{1}{2} \xi \}$.
By (\ref{app89b}) of Lemma \ref{convexupper}, we have $\f{1}{2} \xi
\leq \bs{v}$. This completes the proof of the corollary.

\sect{Proof of Theorem  \ref{lineartheorem} } \la{lineartheoremapp}

To prove Theorem \ref{lineartheorem}, we need a preliminary result.

\beL

\la{lemexist}

Let $\al, \ba, \ze, \eta \in \bb{R}^d$.  If $\bb{E} [ \bs{M} ] \al +
\ze \leq \bb{E} [ S_{\bs{M}} ] \leq \bb{E} [  \bs{M} ] \ba + \eta$,
then there exists $\bs{q} \in \bb{R}^d$ such that $\bs{0}_d \leq
\bs{q} \leq \bs{1}_d$ and that $\bb{E} [ S_{\bs{M}} ] = \bb{E} [
\bs{M} ] \se + \phi$, where $\se = \bs{q} ( \al - \ba ) + \ba $ and
$\phi = \bs{q} (\ze - \eta) + \eta$.

\eeL

\bpf  We consider the scalar case. The argument can be readily generalized to the vector case.
If $\bb{E} [ S_{\bs{M}} ] = \bb{E} [ \bs{M} ] \ba
+ \eta$, then the lemma holds with $\se = \ba$ and $\phi = \eta$.
If $\bb{E} [ S_{\bs{M}} ] = \bb{E} [ \bs{M} ] \al  + \ze$, then the lemma
holds with $\se = \al$ and $\phi = \ze$.
Hence, it remains to prove this lemma under the assumption that \be \la{asplem} \bb{E} [ \bs{M} ] \al +
\ze < \bb{E} [ S_{\bs{M}} ] < \bb{E} [ \bs{M} ] \ba + \eta. \ee For this purpose, define
\[ \se_q = q ( \al - \ba ) + \ba, \qqu  \phi_q = q ( \ze  - \eta ) + \eta, \qqu
w(q) = \bb{E} [ S_{\bs{M}} ] -  \bb{E} [ \bs{M} ] \se_q - \phi_q
\]
for $q \in [0,1]$.  Note that
$w(q) = \bb{E} [ S_{\bs{M}} ] - \bb{E} [ \bs{M} ] [ q ( \al - \ba ) + \ba ] - [ q ( \ze  - \eta ) + \eta ]$ is a
continuous function of $q \in [0, 1]$.  Clearly, as a consequence of (\ref{asplem}),
\[
w(0) = \bb{E} [ S_{\bs{M}} ] - \bb{E} [ \bs{M} ] \ba - \eta < 0,
\qqu w(1) = \bb{E} [ S_{\bs{M}} ] - \bb{E} [ \bs{M} ] \al - \ze > 0.
\]
By virtue of the intermediate value theorem, there exists a number
$q^* \in (0, 1)$ such that $w(q^*) = 0$.  This implies that $\bb{E}
[ S_{\bs{M}} ] =  \bb{E} [ \bs{M} ] \se_{q^*}  + \phi_{q^*}$, where
$\se_{q^*}  = q^* ( \al - \ba ) + \ba$ and $\phi_{q^*} = q^* ( \ze -
\eta ) + \eta$ with $q^* \in (0, 1)$.    This completes the proof of
the lemma.

\epf

\bsk

We are now in a position to prove the theorem.  Since the conditions
(I)--(IV) are fulfilled, it follows from Theorem \ref{BoundGen} that
$\bb{E} [ \bs{M} ] \leq \bb{E} [ \bs{N} ] < \iy$. From
(\ref{boundSM88}), we know that $\bb{E} [ S_{\bs{M}} ]$ exists.
 Since $m > 0$ and $\mscr{R}$ contains $(0,\bs{0}_d)$,  it must be
true that $A s + B t \leq m$ for any $(t, s) \in \mscr{R}$.  Hence,
$\Pr \{  A S_{\bs{M}} + B \bs{M} \leq m \} = 1$.   By Theorem
\ref{fundamental}, we have \be \la{VIPexpression}
 A  \bb{E} [ S_{\bs{M}} ] + B \bb{E} [ \bs{M}
] \leq m. \ee By Lemma \ref{lem1088}, we have $\bb{E} [ \bs{M} ] \al
+ \ze \leq \bb{E} [ S_{\bs{M}} ] \leq \bb{E} [  \bs{M} ] \ba +
\eta$, where
\[
 \al = \mu - \f{\lm}{2} \xi, \qqu \ba =
\mu + \f{\lm}{2} \xi, \qqu \ze = \f{r}{2} \xi, \qqu \eta = -
\f{r}{2} \xi \] with $r = N_0 - K$.  According to Lemma
\ref{lemexist}, there exist $\se^*  = q^* ( \al - \ba ) + \ba$ and
$\phi^* = q^* ( \ze - \eta ) + \eta$, where $q^*$ satisfies
$\bs{0}_d \leq q^* \leq \bs{1}_d$, such that $\bb{E} [ S_{\bs{M}} ]
= \se^* \bb{E} [ \bs{M} ] + \phi^*$. Substituting this expression of
$\bb{E} [ S_{\bs{M}} ]$ into (\ref{VIPexpression}) yields \be
\la{vipineq88}
 \bb{E} [ \bs{M} ] ( A \se^* + B)  + A \phi^* \leq m
\ee Since $(m, m \mu)$ is in the hyperplane $A s + B t = m$, it must
be true that $A\mu + B  = 1$.  As a consequence of the assumption
that $\lm | A | \xi < 2$, we have that \bee A [ \ba + \bs{q} (\al -
\ba) ] + B & = & A \li ( \mu
+ \f{\lm}{2} \xi - \bs{q} \lm \xi \ri )  + B =  1 + \f{\lm}{2} A \xi - \lm A \bs{q} \xi \\
& \geq & 1 + \f{\lm}{2} A \xi - \f{1}{2} \lm (A + | A |)  \xi  =  1
-  \f{1}{2} \lm | A | \xi  > 0 \eee for all $\bs{q}$ such that
$\bs{0}_d \leq \bs{q} \leq \bs{1}_d$. Therefore,  $\min \{ A
\se_{\bs{q}} + B: \bs{0}_d \leq \bs{q} \leq \bs{1}_d \}
> 0$ and  $A \se^* + B > 0$.  Using (\ref{vipineq88}), we have
\[
\bb{E} [ \bs{M} ]  \leq \f{ m - A \phi^* } { B + A \se^* } \leq \max
\li \{ \f{  m - A [ \eta + \bs{q} (\ze - \eta) ]  }{ B +  A [ \ba +
\bs{q} (\al - \ba) ] } : \bs{0}_d \leq \bs{q} \leq \bs{1}_d  \ri \}.
\]
Note that
\[
\f{  m - A [ \eta + \bs{q} (\ze - \eta) ]  }{ B +  A [ \ba + \bs{q}
(\al - \ba) ] } = \f{ m + r A \li ( \f{1}{2} \xi - \bs{q} \xi \ri )
}{  1 + \lm A \li ( \f{1}{2} \xi - \bs{q} \xi \ri )  } = \f{r}{\lm}
+ \f{  m - \f{r}{\lm}  }{  1 + \lm A \li ( \f{1}{2} \xi - \bs{q} \xi
\ri )  }
\]
for all $\bs{q}$ such that $\bs{0}_d \leq \bs{q} \leq \bs{1}_d$.
Recall that under the assumption $\lm | A | \xi < 2$,
\[
1 + \lm A \li ( \f{1}{2} \xi - \bs{q} \xi \ri ) \geq 1 -  \f{1}{2}
\lm | A | \xi  > 0
\]
for all $\bs{q}$ such that $\bs{0}_d \leq \bs{q} \leq \bs{1}_d$.
Hence, {\small $\bb{E} [ \bs{M} ]  \leq \f{r}{\lm} + \f{  m -
\f{r}{\lm} }{ 1 - \f{1}{2} \lm | A | \xi  }$} provided that $m -
\f{r}{\lm} \geq 0$ or equivalently,  $\lm m + K \geq N_0$.  Using
the inequality $\bb{E} [ \bs{N} ] \leq \lm \bb{E} [ \bs{M} ] + K$,
we have \[ \bb{E} [ \bs{N} ] \leq \lm \li [ \f{r}{\lm} + \f{  m -
\f{r}{\lm} }{ 1 - \f{1}{2} \lm | A | \xi   } \ri ] + K = N_0 + \f{
\lm m + K - N_0 }{  1 - \f{1}{2} \lm | A | \xi } \]  provided that
$\lm m + K \geq N_0$.   Of course, we can always set $N_0 = 0$ so
that $\lm m + K \geq N_0$ holds.  This completes the proof of the
theorem.

\sect{Proof of Theorem \ref{ChenLorden} } \la{ChenLordenapp}

Since $A s + B t = m$   is a supporting hyperplane of the continuity
region $\mscr{R}$, passing through the boundary point $(m, m \mu)$,
it must be true that $A \mu + B = 1$. Define
\[
Z = B K + A \sum_{i = 1}^K  X_i \qqu \tx{and} \qqu Z_\ell = B K + A
\sum_{i = (\ell -1) K + 1}^{\ell K}  X_i, \qqu \ell \in \bb{N}.
\]
Then,  $Z_1, Z_2, \cd$ are i.i.d. random vectors having the same
distribution as  $Z$.   By the assumption that each element of
$\bb{E} [ |X|^2]$ is finite, it must be true that each element of
$\bb{E} [ |X|]$ is finite. Hence,
\[
\bb{E} [|Z|]  = \bb{E} \li [ \li | B K + A  \sum_{i=1}^K X_i \ri |
\ri ] \leq |A| \sum_{i=1}^K \bb{E} [ |X_i| ] + |B K| = K ( |A|
\bb{E} [ |X| ] + |B| ),
\]
where the upper bound is finite as a consequence of the assumption
that each element of $\bb{E} [ |X|]$ is finite.  By the definition
of $Z$, we have $\bb{E} [Z] = (A \mu + B) K = K > 0$.   Define
\[
 \bs{\tau} =
\inf \li \{ n \in \bb{N}: \sum_{\ell = 1}^n Z_\ell > m \ri \} \qqu
\tx{and} \qqu \bs{\mcal{N}}  = \inf \{ k \in \mscr{N}: A S_k + B k >
m \},
\] where $\mscr{N} = \{ n \in \bb{N}: \f{n}{K} \in \bb{N} \}$.
Making use of the fact that $\bb{E} [Z] > 0, \; \bb{E} [|Z|] < \iy$
and Theorem \ref{BoundGen}, we have that $\bb{E} [ \bs{\tau} ] <
\iy$. By the definitions of $\bs{\tau}$ and $\bs{\mcal{N}}$, we have
$\bs{\mcal{N}} = K \bs{\tau}$ and thus $\bb{E} [ \bs{\mcal{N}} ] = K
\bb{E} [ \bs{\tau} ] < \iy$. Since $m
> 0$ and the convex set $\mscr{R}$ contains $(0, \bs{0}_d)$, it must
be true that $A s + B t \leq m$ for all $(t, s) \in \mscr{R}$.
Hence, $\{ (n, S_n) \in  \mscr{R} \} \subseteq \{   A S_n + B n \leq
m \}$ for $n \in \bb{N}$.  This implies that $\bs{N} \leq
\bs{\mcal{N}}$. Hence, $\bb{E} [ \bs{N} ] \leq \bb{E} [
\bs{\mcal{N}} ] < \iy$.  Note that \bel \bb{E} [Z^2]  & = & \bb{E}
\li [ \li | B K + A \sum_{i=1}^K X_i \ri |^2 \ri ]  =
 \bb{E} \li [ \li |A \sum_{i=1}^K  X_i \ri |^2 \ri ] + 2 B K \bb{E} \li [ A  \sum_{i=1}^K X_i \ri ] + (B K)^2 \nonumber\\
& = &   \bb{E} \li [ \li | A \sum_{i=1}^K X_i \ri |^2 \ri ] + 2 B K^2  A \mu + (B K)^2 \nonumber\\
& \leq &  ||A||^2 \times  \sum_{i=1}^K \bb{E} [|| X_i ||^2] + 2 B
K^2 A \mu + (B K)^2 \la{useCW888}\\
& = & K  ||A||^2 \times \bb{E} [ || X ||^2 ]  + 2 B K^2 A \mu + (B
K)^2, \la{z2bound} \eel where we have used Cauchy-Schwarz inequality
to obtain (\ref{useCW888}).   By the assumption that each element of
$\bb{E} [ | X |^2 ]$ is finite, it must be true that $\bb{E} [ || X
||^2 ]$ is finite.   Using the bound in (\ref{z2bound}) and the
finiteness of  $\bb{E} [ || X ||^2 ]$, we have that $\bb{E} [Z^2] <
\iy$.  By Lorden's inequality
\[
\bb{E} \li [ \sum_{\ell = 1}^{\bs{\tau}}  Z_\ell -  m \ri ] \leq \f{
\bb{E} [ (Z^+)^2 ] }{ \bb{E} [ Z ]  }.
\]
Using Wald's equation, we have $\bb{E} \li [ \sum_{\ell =
1}^{\bs{\tau}}  Z_\ell \ri ] = \bb{E}[ \bs{\tau} ] \bb{E} [ Z ]$ and
thus $\bb{E}[ \bs{\tau} ] \bb{E} [ Z ]  - m \leq \f{ \bb{E} [
(Z^+)^2 ] }{ \bb{E} [ Z ]  }$,  from which we have
\[
\bb{E}[ \bs{\tau} ]  \leq \f{ m }{ \bb{E} [ Z ] } + \f{ \bb{E} [
(Z^+)^2 ] }{ \bb{E}^2 [ Z ]  }.
\]
Hence, \bel
 \bb{E}[ \bs{N} ]  & \leq & \bb{E}[ \bs{\mcal{N}} ] = K \bb{E}[ \bs{\tau} ] \leq
  \f{ m K }{ \bb{E} [ Z ] } + \f{K \bb{E} [ (Z^+)^2 ] }{ \bb{E}^2 [ Z ]  } \nonumber\\
 & = & m +
\f{ \bb{E} [ (Z^+)^2 ] }{ K  }, \la{uselater} \eel where we have
used the fact that $\bb{E} [ Z ] = K$.  Note that \bel \bb{E} [ (Z^+)^2 ] & \leq & \bb{E} [ Z^2 ]
 =  \bb{E} \li [ \li ( B K + \sum_{i = 1}^K A X_i \ri )^2 \ri ]
  =  \bb{E} \li [ \li ( B K + K A \mu + \sum_{i = 1}^K A ( X_i - \mu) \ri )^2 \ri ] \nonumber\\
& = & K^2 ( A \mu + B )^2 +  \bb{E} \li [ \li (  \sum_{i = 1}^K A ( X_i - \mu) \ri )^2 \ri ] \nonumber \\
& = & K^2  +  K  \bb{E} \li [ \li |  A ( X - \mu) \ri |^2 \ri ],
\la{usenow8899} \eel where we have used the assumption that $X_1,
X_2, \cd$ are identically distributed and mutually independent
random vectors.   It follows from (\ref{uselater}) and
(\ref{usenow8899}) that  \be \la{usenow8868}
 \bb{E}[ \bs{N} ]  \leq
m + K + \bb{E} \li [ \li | A ( X - \mu) \ri |^2 \ri ]. \ee  Using
Cauchy-Schwarz inequality, we have
\[
\bb{E} \li [ \li | A ( X - \mu) \ri |^2 \ri ] \leq \bb{E} \li [  ||
A ||^2 \times || X - \mu ||^2  \ri ] =  || A ||^2 \times \bb{E} \li
[ || X - \mu ||^2  \ri ].
\]
Therefore, $\bb{E}[ \bs{N} ]  \leq m + K + \bb{E} \li [ \li | A ( X
- \mu) \ri |^2 \ri ] \leq m + K +  || A ||^2 \times \bb{E} \li [ ||
X - \mu ||^2 \ri ]$.  This establishes assertion (I) of the theorem.

If  the  elements of $X$ are mutually independent, then $\bb{E} \li
[ \li | A ( X - \mu) \ri |^2 \ri ] = A^2 \bb{E} \li [ | X - \mu |^2
\ri ]$. It follows from this fact and (\ref{usenow8868}) that
$\bb{E}[ \bs{N} ]  \leq m + K +  A^2 \; \bb{E} [ | X - \mu |^2 ]$.
This establishes assertion (II) of the theorem.

It remains to show assertion (III).  As a consequence of the
definition of $Z$ and the assumption that $\Pr \{ \bs{a} \leq X \leq
\bs{b} \} = 1$, we have $K u \leq Z \leq K v$ almost surely.  It
follows that \bee (Z^+)^2 \leq Z^2 \leq  \f{ (Kv)^2 - (K u)^2 }{ K v
- K u } (Z - K u) + (K u)^2 = K (u + v) Z - K^2 u v \eee  almost
surely.  Hence, \be \la{usebefore}
 \bb{E} \li [  (Z^+)^2 \ri ] \leq K (u + v) \bb{E} [ Z ]  - K^2 u v = K^2 (u + v - u v).
\ee Making use of (\ref{uselater}) and (\ref{usebefore}), we have
\bee \bb{E} [ \bs{N} ]  \leq  m + \f{ \bb{E} \li [  (Z^+)^2 \ri ] }{
K  }  =   m + \f{ K^2 (u + v  - u v )}{ K  }  =   m + K (u + v - u v
).  \eee  This establishes the first inequality of assertion (III).

To show the second inequality of assertion (III), note that
\[
(Z^+)^2 \leq \f{ (K v)^2 }{ K v  - K u } (Z - K u)  = \f{ K v^2}{v -
u}  (Z - K u)
\]
almost surely for $u < 0$.  Hence, \be \la{usebeforeb}  \bb{E} \li [
(Z^+)^2 \ri ] \leq \f{ K v^2}{v - u}  (\bb{E} [ Z ] - K u) = \f{K^2
v^2}{v - u} ( 1 - u), \qqu u < 0. \ee Making use of (\ref{uselater})
and (\ref{usebeforeb}), we have \bee \bb{E} [ \bs{N} ]  \leq  m +
\f{ \bb{E} \li [  (Z^+)^2 \ri ] }{ K  }  =   m + \f{  \f{K^2 v^2}{v
- u}  ( 1 - u) }{ K  }  =   m +  K v^2  \li ( \f{ v - u }{ 1 - u }
\ri ) \eee for $u < 0$.

This completes the proof of the theorem.

\sect{Proof of Theorem \ref{discrerebound} } \la{discrereboundap}

We need a preliminary result.

\beL

\la{chenconvex}

  Let
$(m, m \mu)$, where $m = \mscr{B} (\mu)
> 0$, be a boundary point of the continuity region $\mscr{R}$.
Let $\nabla (v)$ denote the gradient of $\ln \mscr{B} (v)$. Define
\[
A = - \nabla (\mu), \qqu B = 1 - A \mu, \qqu C = m.
\]
Then,  $A s  + B t = C$ is the supporting hyperplane for $\mscr{R}$
passing through the boundary point $ (m,  m \mu)$. \eeL

\bpf

Define function
\[
f (t, s) = t - \mscr{B} \li ( \f{s}{t} \ri )
\]
for $(t, s) \in \mscr{R}$ with $t > 0$.  As a consequence of  the
convexity of $\mscr{R}$ and the definition of the function $\mscr{B}
(.)$, it must be true that $f (t, s) = 0$ holds for any boundary
point $(t, s)$ of $\mscr{R}$ with $t > 0$. In particular, $f(m, m
\mu) = 0$. Since $\mscr{B} (v)$ is differentiable at  $v = \mu$, the
function $f (t, s)$ is differentiable at $(t, s) = (m,  m u)$. Since
$\mscr{R}$ is convex, it must be true that the tangent plane to the
surface $f(t,s) = 0$, passing through $(m,m \mu)$, coincides with
the supporting hyperplane of $\mscr{R}$. Therefore, to show that $A
s + B t = C$ is the supporting hyperplane of $\mscr{R}$ passing
through the boundary point $ (m, m \mu)$, it suffices to show that
$C$ is equal to $m$, and that $A$ and $B$ are, respectively, equal
to the partial derivatives of $f (t, s)$ with respect to $s$ and $t$
at $(t, s) = (m, m \mu)$. In other words, it is sufficient to show
that
\[
C = m, \qqu A = \li. \f{ \pa f (t, s) } {\pa s} \ri |_{ t = m, \; s
= m \mu  }, \qqu B = \li. \f{ \pa f (t, s) } {\pa t} \ri |_{ t = m,
\; s = m \mu  }.
\]
Define $h (v) = \f{ \pa \mscr{B} (v) } {\pa v}$.  Using the chain
rule of differentiation, we have
\[
\f{ \pa f (t, s) } {\pa s} = - \f{1}{t} h \li (\f{s}{t} \ri) \qqu
\tx{and} \qqu \f{ \pa f (t, s) } {\pa t} = 1  +   h \li ( \f{s}{t}
\ri ) \f{s}{t^2}.
\]
Evaluating such derivatives with $t = m, \; s = m \mu$ yields
\[
\li. \f{ \pa f (t, s) } {\pa s} \ri |_{ t = m, \; s = m \mu  }= -
\f{ h (\mu) }{ \mscr{B} (\mu)} = - \nabla  (\mu ) = A
\]
and
\[
\li. \f{ \pa f (t, s) } {\pa t} \ri |_{ t = m, \; s = m \mu  }= 1 +
 \f{ h (\mu) }{ \mscr{B} (\mu)} \mu  = 1 - A \mu  = B.
\]
Since the boundary point $ (m,  m \mu)$ is in the supporting
hyperplane, it must be true that
\[
C = m A \mu + B m = m (   A \mu + B ).
\]
Observing that $A \mu + B = 1$, we have $C = m$. This completes the
proof of the lemma.

\epf

We are now in a position to prove the theorem. Making use of Theorem
\ref{ChenLorden} and Lemma \ref{chenconvex}, we have
$$\bb{E}[ \bs{N} ]  \leq m + K +  \bb{E} \li [ \li | A (X - \mu) \ri |^2 \ri ]
\leq m + K +  || A ||^2 \times \bb{E} [ || X - \mu
||^2],$$ where $A = - \nabla  (\mu) = - V$ and $m = \mscr{B} (\mu)$.
It follows that
$$\bb{E}[ \bs{N} ]  \leq  \mscr{B} (\mu) + K + \bb{E} \li [ \li | V (X - \mu) \ri |^2 \ri ]
\leq \mscr{B} (\mu) + K + || V ||^2 \times \bb{E} [ || X - \mu
||^2].$$ This completes the proof of the theorem.

\sect{Proof of Corollary \ref{concavefunction} }
\la{concavefunctionapp}

Note that the stopping time $\bs{N}$ defined by (\ref{exam8888}) can
be expressed in the more general form (\ref{newST8896}) with
continuity region $\mscr{R} = \{ (t, s): t \in \bb{R}^+, \; s \in
\bb{R}, \; s \leq f (t) \}$.   Let the DET function associated with
$\mscr{R}$ be $\mscr{B} (v)$. Then, the solution of the equation $m
\mu = f (m)$ can be taken as $m = \mscr{B} (\mu)$. Since $f(t)$ is
differentiable at $t = m$, it follows that $\mscr{B} (v)$ must be
differentiable at $v = \mu$. Due to the concavity of the boundary
function $f(.)$, the continuity region $\mscr{R}$ is a convex set.
To apply Theorem \ref{discrerebound} to bound the
 stopping time in (\ref{exam8888}), we can calculate  $\nabla(\mu)$ by (\ref{defnabla}) as follows.

By the definition of the DET function $\mscr{B} (v)$,  we have
$\mscr{B} (v) v = f ( \mscr{B} (v) )$ at a neighborhood of $v =
\mu$. Differentiating both sides of this equation with respect to
$v$ at $v = \mu$  yields
\[
\mscr{B}^\prime (\mu) \mu + \mscr{B} (\mu) = f^\prime (\mscr{B}
(\mu) ) \mscr{B}^\prime (\mu),
\]
where $f^\prime (.)$ and $\mscr{B}^\prime (.)$ denote the first
derivatives of $f (.)$ and $\mscr{B}(.)$, respectively.   It follows
that the first derivative, $\mscr{B}^\prime (\mu)$, of $\mscr{B}
(v)$ at $v = \mu$ can be obtained as
\[
\mscr{B}^\prime (\mu) = \f{ \mscr{B} (\mu)  }{  f^\prime (\mscr{B}
(\mu) ) - \mu } = \f{m}{ f^\prime (m) - \mu }.
\]
Therefore, the gradient of $\ln \mscr{B} (v)$ at $v = \mu$ is
\[
\nabla(\mu) = \f{ \mscr{B}^\prime (\mu)  } {  \mscr{B} (\mu) } =
\f{m}{ [f^\prime (m) - \mu ] \mscr{B} (\mu) } = \f{1}{ f^\prime (m)
- \mu }.
\]
It follows from Theorem \ref{discrerebound} that $\bb{E} [ \bs{N} ]
\leq m + K + \f{ \si^2  }{  | f^\prime (m) - \mu |^2 }$.

\sect{Proof of Theorem \ref{BoundConcentrationa88} }
\la{BoundConcentrationa88app}

 Note that
\bel
 \bb{E} [ \bs{N} ] & = & \sum_{i=0}^\iy \Pr \{ \bs{N} > i \}
=  N_1 + \sum_{\ell = 1}^\iy  (N_{\ell + 1} - N_\ell ) \Pr \li \{
\bs{N} >
N_\ell \ri \} \nonumber\\
& = & N_1 + \sum_{\ell = 1}^\ka  (N_{\ell + 1} - N_\ell ) \Pr \li \{
\bs{N} > N_\ell \ri \}.  \la{VIP8899a} \eel By the definition of
$\bs{N}$ and $\de (t)$, we have \bee \li \{  \bs{N} > N_\ell \ri \}
\subseteq  \{ (N_\ell, S_{N_\ell}) \in \mscr{R}  \} \subseteq  \li
\{ || \ovl{X}_{N_\ell} - \mu || \geq \de (N_\ell) \ri \} \eee for
$\ell \in \bb{N}$.  It follows  that \be \Pr \li \{  \bs{N} > N_\ell
\ri \} \leq  \Pr \li \{ || \ovl{X}_{N_\ell} - \mu || \geq \de
(N_\ell) \ri \} \la{VIP8899b} \ee for $\ell \in \bb{N}$.  Making use
of (\ref{VIP8899a}) and (\ref{VIP8899b}), we have $\bb{E} [ \bs{N} ]
\leq N_1 + \sum_{\ell = 1 }^\ka  (N_{\ell + 1} - N_\ell )
 \Pr \li \{  \li | \ovl{X}_{N_\ell} -  \mu \ri | \geq \de (N_\ell) \ri
 \}$.    This completes the proof of the theorem.

\sect{Proof of Theorem  \ref{BoundConcentrationb88} }
\la{BoundConcentrationb88app}

We need a preliminary result.

\beL

\la{min8899}

Let $A$ be a row matrix of size $1 \times d$ and $B, \; C \in
\bb{R}$. For $t \in \bb{R}$ and $\mu \in \bb{R}^d$,
\[
\min \{ || s - t \mu ||:  s \in \bb{R}^d, \;  A s + B
 t = C \} = \f{ | C - t (A \mu + B) | }{
|| A || }
\]
provided that $|| A || > 0$.

\eeL

\bpf

Note that $A (s - t \mu) = C - t ( A \mu + B )$ holds for any point
$(t, s)$ in the hyperplane $A s + B t = C$.  Hence, \[ | C - t ( A
\mu + B ) |  = | A (s - t \mu) | \leq || A || \times || s - t \mu
||,
\]
where the inequality follows from Cauchy-Schwarz inequality.  Since
$||A|| > 0$, we have $|| s - t \mu || \geq \f{ | C - t (A \mu + B) |
}{ || A || }$ for any point $(t, s)$ in the hyperplane. Now let $s^*
= t \mu + \f{ [ C - t (A \mu + B) ] A^\top }{ || A ||^2 }$. It can
be checked that \[ || s^* - t u || = \f{ | C - t (A \mu + B) | }{ ||
A || }
\]
and that $A s^* + B t = C$, which implies that the  minimum of $|| s
- t \mu ||$ is attained at point $(t, s^*)$ of the hyperplane.  This
completes the proof of the lemma.

\epf

 We are now in a position to prove the theorem.  It is sufficient to consider two
 cases as follows.

 Case(i): $|| A || = 0$.

 Case (ii): $|| A || > 0$.

In Case (i), we have $B = 1$, since the hyperplane contains point
$(m, m \mu)$. Hence, $\mscr{R} \subseteq \{ (t, s): 0 \leq t \leq m,
\; s \in \bb{R}^d  \}$. It follows that $\bs{N} = N_\jmath$. On the
other hand, the right side of the inequality
(\ref{inelinear8899883}) is equal to $N_\jmath$ because {\small $\Pr
\li \{ || A || \times \li | \li | \ovl{X}_{N_\ell} -  \mu \ri | \ri
| \geq 1 - \f{m}{N_\ell}  \ri \} = 0$} for all $\ell \geq \jmath$.
Hence, Theorem
 \ref{BoundConcentrationb88} holds trivially in Case (i). It remains
 to show the theorem in Case(ii). We proceed as follows.

 By Theorem
 \ref{BoundConcentrationa88},
\bel \bb{E} [ \bs{N} ] & \leq & N_1 + \sum_{\ell = 1 }^\ka  (N_{\ell
+ 1} - N_\ell )  \Pr \li \{  \li | \ovl{X}_{N_\ell} -  \mu \ri |
\geq \de (N_\ell) \ri  \} \nonumber\\
& \leq & N_\jmath + \sum_{\ell = \jmath }^\ka  (N_{\ell + 1} -
N_\ell ) \Pr \li \{  \li | \ovl{X}_{N_\ell} -  \mu \ri | \geq \de
(N_\ell) \ri \}. \la{last668899}
 \eel
For $t > 0$, \bel  \de (t) & = & \inf \{ || v - \mu ||: v \in
\bb{R}^d, \; (t , v t) \in \mscr{R} \} =  \f{1}{t} \times \inf \{ ||
s - t \mu ||: s \in \bb{R}^d, \; (t , s)
\in \mscr{R} \} \nonumber \\
& \geq & \f{1}{t} \times \inf \{ || s - t \mu ||: s \in \bb{R}^d, \;
A s + B t \leq m  \}, \la{explain}  \eel   where (\ref{explain}) is
due to the fact that $\mscr{R} \subseteq \{ (t, s): t \geq 0, \; s
\in \bb{R}^d, \; As + B t \leq m \}$.   Since the supporting
hyperplane $A s + B t = m$ contains $(m, m \mu)$, we have $A \mu + B
= 1$.  For $t \geq N_\jmath$ and $s = \mu t$, we have $As + B t = t
( A \mu + B ) = t \geq N_\jmath
> m$. This implies that $(t, t \mu) \in \{ (t, s): s \in \bb{R}^d,
\; As + B t
> m \}$ for $t \geq N_\jmath$. For any point $(t, z) \in \{ (t, s): s
\in \bb{R}^d, \; As + B t <  m \}$,  there exists a unique point,
$(t, s^*)$, of the hyperplane $As + B t = m$ such that $s^* = r (t
\mu) + (1 - r)  z$ for some number $r \in (0,1)$.  Hence, $|| s^* -
\mu t || = || r (t \mu) + (1 - r)  z - \mu t || = (1 - r ) || z -
\mu t|| \leq ||z - \mu t||$.  It follows that \be \la{key663689}
 \inf \{ || s - t \mu ||: s \in \bb{R}^d, \; A s + B t
\leq m \} = \inf \{ || s - t \mu ||: s \in \bb{R}^d, \; A s + B t =
m \} \ee for $t \geq N_\jmath$.   According to Lemma \ref{min8899},
we have \be \la{8899368}
 \inf \{ || s - t \mu ||: s \in \bb{R}^d, \; A
s + B t = m  \} = \f{ | m - t (A \mu + B) | }{ || A || } = \f{ | m -
t | }{ || A || }. \ee  Making use of (\ref{explain}),
(\ref{key663689}) and (\ref{8899368}), we have $\de (t) \geq
\f{1}{t} \f{ | m - t | }{ || A || } = \li ( 1 - \f{m}{t} \ri )
\f{1}{ ||A|| }$ for $t \geq N_\jmath$.  Therefore,  \be
\la{usenow88}
 \de (N_\ell) \geq  \li ( 1 - \f{m}{N_\ell} \ri ) \f{1}{
||A|| }. \ee for $\ell \geq \jmath$. Making use of
(\ref{last668899}) and (\ref{usenow88})  completes the proof of the
theorem.

\sect{Proof of Theorem   \ref{scalarAFPT}} \la{scalarAFPTapp}

By Theorem \ref{BoundConcentrationa88},  \be \la{scalineq8899}
 \bb{E} [ \bs{N} ]  \leq   N_1 + \sum_{\ell = 1}^\ka  (N_{\ell + 1} - N_\ell ) \Pr \li \{
\bs{N} > N_\ell \ri \}  \leq   N_\jmath + \sum_{\ell = \jmath}^\ka
(N_{\ell + 1} - N_\ell ) \Pr \li \{ \bs{N} > N_\ell \ri \}.  \ee By
the definition of $\bs{N}$, we have $\li \{ \bs{N}
> N_\ell \ri \} \subseteq \{ (N_\ell, S_{N_\ell}) \in \mscr{R}  \}$ for $\ell \in \bb{N}$.
Since $X$ is a scalar random variable and the continuity region
$\mscr{R}$ is convex, it follows from the definition of $\de(t)$
that \bel
 &  &  \{ (N_\ell, S_{N_\ell}) \in
\mscr{R}  \} \subseteq \{ \ovl{X}_{N_\ell} \geq \mu + \de (N_\ell)
\} \qqu  \tx{if} \qu \{ (t, s) \in \mscr{R}: t >
m, \; s > \mu t \} \neq \emptyset, \la{twocasesa}\\
&  &  \{ (N_\ell, S_{N_\ell}) \in \mscr{R}  \} \subseteq  \{
\ovl{X}_{N_\ell} \leq \mu - \de (N_\ell) \} \qqu  \tx{if} \qu \{ (t,
s) \in \mscr{R}: t > m, \; s < \mu t \} \neq \emptyset
\la{twocasesb} \eel for all $\ell \geq \jmath$.  Making use of
(\ref{scalineq8899}),  (\ref{twocasesa}) and (\ref{twocasesb})
proves the inequalities (\ref{ineq2879a}) and (\ref{ineq2879b}) of
the theorem.  It remains to show the inequality (\ref{ineq2882})
under additional assumption that there exists a supporting
hyperplane $A s + B t = m$ of $\mscr{R}$,  passing through $(m, m
\mu)$.  For this purpose, it suffices to consider three cases.

Case (i): $A = 0$.

Case (ii): $A > 0$.

Case (iii): $A < 0$.

In Case (i), we have $B = 1$, since the hyperplane contains point
$(m, m \mu)$. Hence, $\mscr{R} \subseteq \{ (t, s): 0 \leq t \leq m,
\; s \in \bb{R}  \}$. It follows that $\bs{N} = N_\jmath$. On the
other hand, the right side of the inequality (\ref{ineq2882}) is
equal to $N_\jmath$ because {\small $\Pr \li \{  A \li (
\ovl{X}_{N_\ell} - \mu \ri )  \geq  1 - \f{m}{N_\ell}   \ri \}  =
0$} for all $\ell \geq \jmath$. Hence, the inequality
(\ref{ineq2882}) holds trivially in Case (i).

In Case (ii),  $|A| = A > 0$.  Using the same argument as that in
the proof of Theorem \ref{BoundConcentrationb88}, we have {\small
$\de (N_\ell) = \f{1}{A} \li ( 1 - \f{m}{N_\ell}  \ri )$} for $\ell
\geq \jmath$. Hence, {\small $\Pr \{ \ovl{X}_{N_\ell} \geq \mu + \de
(N_\ell) \} = \Pr \li \{  A \li ( \ovl{X}_{N_\ell} -  \mu \ri ) \geq
1 - \f{m}{N_\ell}   \ri \}$} for all $\ell \geq \jmath$.  Making use
of this fact and (\ref{ineq2879a}) shows (\ref{ineq2882}) for Case
(ii).

In Case (iii),  $|A| = - A > 0$.  Using the same argument as that in
the proof of Theorem \ref{BoundConcentrationb88}, we have {\small
$\de (N_\ell) = - \f{1}{A} \li ( 1 - \f{m}{N_\ell}  \ri )$} for
$\ell \geq \jmath$. Hence, {\small $\Pr \{ \ovl{X}_{N_\ell} \leq \mu
- \de (N_\ell) \}  = \Pr \li \{  A \li ( \ovl{X}_{N_\ell} - \mu \ri
) \geq  1 - \f{m}{N_\ell}   \ri \}$} for all $\ell \geq \jmath$.
Making use of this fact and (\ref{ineq2879b}) shows (\ref{ineq2882})
for Case (iii).

This completes the proof of the theorem.

\sect{Proof of Theorem \ref{Brown2} } \la{Brown2app}

We shall first show $\bb{E} [ T ] \geq \mscr{A} (\mu)$ under the
assumption that $\mscr{A} (\mu ) < \iy$. If $\bb{E} [ T ] = \iy$,
then $\bb{E} [ T ] \geq \mscr{A} (\mu)$ trivially holds. If $\bb{E}
[ T ] < \iy$, then $\Pr \{ T < \iy \} = 1$ and it follows that
$X_{T}$ is well-defined and $\Pr \{ (T, X_{T}) \in \mscr{R}^c \} =
1$.  Since $\bb{E} [ T ] < \iy$, it follows from Wald's equation
that $\bb{E} [  X_{T}] = \bb{E} [ T ] \mu$.  According to Theorem
\ref{fundamental}, we have $( \bb{E} [ T ], \;  \bb{E} [  X_{T}] )
\in \mscr{R}^c$.   Hence, $( \bb{E} [ T ], \;  \bb{E} [ T ] \mu  )
\in \mscr{R}^c$. It follows from the notion of IDET that $\bb{E} [ T
] \geq \mscr{A} (\mu)$.  This establishes the first assertion.

It remains to show that $\bb{E} [ T ] = \iy$ under the assumption
that $\mscr{A} (\mu) = \iy$.  We use a contradiction method. Suppose
that $\bb{E} [ T ] < \iy$, then $\Pr \{ T < \iy  \} = 1$ and it
follows that $\Pr \{ (T, X_{T}) \in \mscr{R}^c \} = 1$. Since
$\bb{E} [ T ] < \iy$, it follows from Wald's equation that $\bb{E} [
X_{T}] = \bb{E} [ T ] \mu$.  According to Theorem \ref{fundamental},
we have $( \bb{E} [ T ], \;  \bb{E} [ X_{T}] ) \in \mscr{R}^c$.
Hence, $( \bb{E} [ T ], \;  \bb{E} [ T ] \mu  ) \in \mscr{R}^c$,
which immediately implies that $\mscr{A}(\mu) \leq \bb{E} [ T ] <
\iy$. This is a contradiction. Therefore, it must be true that
$\bb{E} [ T ] = \iy$ if $\mscr{A} (\mu) = \iy$.  The proof of the
theorem is thus completed.

\sect{Proof of Corollary \ref{Brown4} } \la{Brown4app}

 Note that $T = \inf \{ t > 0: (t, X_t) \notin
\mscr{R} \}$, where $\mscr{R} = \{ (t, s): t > 0, \; s \in \bb{R}^d,
\; t g (\f{s}{t}) \leq 1 \}$.  Since $g$ is a concave function on
$\bb{R}^d$, it follows from Lemma \ref{BoundaryCon} that $t g
(\f{s}{t})$ is a concave function of $t > 0, \; s \in \bb{R}^d$.
Hence, $\mscr{R}^c = \{ (t, s): t > 0, \; s \in \bb{R}^d, \; t g
(\f{s}{t}) > 1 \}$ is a convex set.   Note that
\[
\mscr{A} (\mu) = \inf \{ t > 0: (t, \mu t) \notin \mscr{R} \} = \inf
\{ t > 0: (t, \mu t) \in \mscr{R}^c \} = \inf \{ t > 0: t g (\mu) >
1 \} = \f{1}{g (\mu)} < \iy.
\]
It follows from Theorem \ref{Brown2} that $\bb{E} [ T ] \geq \f{1}{g
(\mu)}$.  This completes the proof of the corollary.

\section{Proof of Theorem \ref{them8899l}} \la{them8899lpf}

We need a preliminary result.

\beL

\la{convergeas}

Let $\mcal{T} = \inf \{ t > 0: (t, X_t) \notin \mcal{R} \}$, where
$\mcal{R}$ is a closed  subset of $\{ (t, s): t \in \bb{R}^+, \; s
\in \bb{R}^d \}$ which contains $(0, \bs{0}_d)$. Define a sequence
of random variables $\{ \mcal{T}_k, \; k \in \bb{N} \}$ such that
$\mcal{T}_k = \inf \{ t > 0:  (t, X_t) \notin \mcal{R}, \;
\tx{where} \; t = i 2^{- k} \; \tx{with} \; i \in \bb{N} \}$ for $k
\in \bb{N}$. Then, $\mcal{T}_k \to \mcal{T}$ almost surely as $k \to
\iy$.

\eeL

\bpf

Since $\{ X_t, \; t \geq 0 \}$ is a L\'{e}vy process, there exists
$\Om_0 \in \Om$ such that $\Pr \{ \Om_ 0 \} = 1$ and that for every
$\om \in \Om_0$, the sample path $X_t(\om)$ is right-continuous for
all $t > 0$. Hence, to show the lemma, it suffices to show that for
every $\om \in \Om_0$, $\mcal{T}_k (\om) \to \mcal{T} (\om)$ as $k
\to \iy$.

Since $\mcal{R}$ is closed, it follows that the complementary set,
$\mcal{R}^c$,  of $\mcal{R}$ must be open. Let $\om \in \Om_0$ and
$\varsigma = \mcal{T} (\om), \; x_\varsigma = X_{\mcal{T}} (\om)$.
Then, $(\varsigma, x_\varsigma ) \in \mcal{R}^c$. Since $\mcal{R}^c$
is an open set, it follows that there exists $\eta > 0$ such that
$\{ (t, z): (t - \varsigma)^2 + || z - x_\varsigma ||^2 < \eta^2 \}
\subset \mcal{R}^c$. By the right-continuity of the sample paths of
a L\'{e}vy process, there exists $\vep \in (0, \f{\eta}{\sq{2}})$
such that $|| X_{\varsigma + \de} (\om) - x_\varsigma || <
\f{\eta}{\sq{2}}$ for $0 < \de < \vep$.  Hence,
\[
(\varsigma + \de - \varsigma)^2 + || X_{\varsigma + \de} (\om) -
x_\varsigma ||^2 < \li ( \f{\eta}{\sq{2}} \ri )^2 + \li (
\f{\eta}{\sq{2}} \ri )^2 = \eta^2
\]
for $0 < \de < \vep$.   This implies that \be \la{usestar}
 (\varsigma +
\de, X_{\varsigma + \de} (\om) ) \in \{ (t, z): (t - \varsigma)^2 +
|| z - x_\varsigma ||^2  < \eta^2 \} \subset \mcal{R}^c \qqu
\tx{for} \qu 0 < \de < \vep. \ee By the definition of $\mcal{T}_k$,
we have that $\mcal{T}_k (\om) = \inf \{ t \geq \varsigma: (t, X_t)
\notin \mcal{R}, \; \tx{where} \; t = i 2^{- k} \; \tx{with} \; i
\in \bb{N} \}$ for $k \in \bb{N}$, where $\de_k = 2^{-k}$. Clearly,
$\de_k < \vep$ for $k > \log_2 \f{1}{\vep}$. Therefore, for $k >
\log_2 \f{1}{\vep}$, it follows from (\ref{usestar}) that $(
\varsigma + \de_k, X_{\varsigma + \de_k} ) \notin \mcal{R}$,  which
implies that $\varsigma \leq \mcal{T}_k (\om) \leq \varsigma +
\de_k$ for all $k
> \log_2 \f{1}{\vep}$.  It follows that $\mcal{T}_k (\om) \to \varsigma$
as $k \to \iy$.  This proves that $\mcal{T}_k \to \mcal{T}$ almost
surely as $k \to \iy$.

\epf

 We are now in a position to prove the theorem.  Let $\de_k = 2^{-k}$ for $k \in \bb{N}$.
 Define $Y_n = X_{n \de_k} - X_{(n - 1) \de_k}$ and $Z_n = \sum_{i=1}^n Y_i$ for $n \in \bb{N}$.
As a consequence of the stationary independent increments property
of a L\'{e}vy process, $\{ Y_n, \; n \in \bb{N} \}$ is a sequence of
i.i.d. random vectors of common mean $\mu \de_k$.    Define
\[
\mcal{T}_k = \inf \{ t > 0:  A X_t + B t > \tau,  \; \tx{where} \; t
= n \de_k \; \tx{with} \;  n \in \bb{N} \}, \qqu k \in \bb{N}
\]
and $\mcal{M}_k = \inf  \{ n \in \bb{N}: (n, Z_n) \notin \mcal{R}
\}$ for $k \in \bb{N}$,  where
\[
\mcal{R} = \li \{ (t, s): t \geq 0, \; s \in \bb{R}^d, \;
\f{1}{\de_k} A s + B t \leq \f{\tau}{\de_k} \ri \}, \qqu k \in
\bb{N}.
\]
It can be checked that $\mcal{T}_k = \de_k \times  \mcal{M}_k$ for
$k \in \bb{N}$.   Note that $\mcal{R}$ is a convex continuity region
associated with $\mcal{M}_k$, with DET function
\[
g (v) = \sup \{ t > 0: (t , v t) \in \mcal{R} \} = \f{\tau}{
 A v + \de_k B  }, \qqu v \in \bb{R}^d.
\]
Since $A s + B t = \tau$ contains point $(\tau, \mu \tau)$, it must
be true that $A \mu + B = 1$. Hence, $g (\mu \de_k) =
\f{\tau}{\de_k}$.   Note that $\f{1}{\de_k} A s + B t =
\f{\tau}{\de_k}$ is a supporting hyperplane of $\mcal{R}$,  passing
through point $ ( g (\mu \de_k), g (\mu \de_k) \mu \de_k ) = \li (
\f{\tau}{\de_k},  \tau \mu \ri )$.  Making use of assertion (I) of
Theorem \ref{ChenLorden}, we have
\[
\bb{E}[  \mcal{M}_k ]  \leq  \f{\tau}{\de_k} + 1 +  \bb{E} \li [ \li
| \f{A}{\de_k}  ( X_{\de_k}  - \mu \de_k) \ri |^2 \ri ] \leq
\f{\tau}{\de_k} + 1 + \li | \li | \f{A}{\de_k} \ri | \ri |^2 \times
\bb{E} \li [ || X_{\de_k} - \mu \de_k ||^2 \ri ].
\]
Hence, \[ \bb{E} [ \mcal{M}_k ]  \leq  \f{\tau}{\de_k} + 1 + \f{1}{
(\de_k )^2 } \bb{E} \li [  \li | A ( X_{\de_k} - \mu \de_k ) \ri |^2
\ri ] \leq \f{\tau}{\de_k}   + 1 + \f{1}{ ( \de_k )^2 } \li | \li |
A \ri | \ri |^2 \times \bb{E} [ || X_{\de_k} - \mu \de_k ||^2 ].
\]
Making use of the above inequalities and the relation $\mcal{T}_k =
\mcal{M}_k \de_k$, we have \be \la{use88993a}
 \bb{E} [ \mcal{T}_k ] \leq
\tau + \de_k + \f{1}{\de_k} \bb{E} \li [  \li | A ( X_{\de_k} - \mu
\de_k ) \ri |^2 \ri ] \leq \tau  + \de_k + \f{1}{ \de_k } \li | \li
| A \ri | \ri |^2 \times \bb{E} [ || X_{\de_k} - \mu \de_k ||^2 ].
\ee Since $\{ X_t, \; t \geq 0 \}$ is a L\'{e}vy process, we have
\be \la{use88993b}
 \bb{E} [ || X_{\de_k} - \mu \de_k  ||^2 ] = \de_k \bb{E} [ || X - \mu ||^2 ]. \ee Since $\{
X_t, \; t \geq 0 \}$ is a L\'{e}vy process,  it follows that $\{ A (
X_t - \mu t), \; t \geq 0 \}$ is also a L\'{e}vy process and thus
\be \la{use88993c}
 \bb{E} \li [  \li | A ( X_{\de_k} - \mu \de_k ) \ri
|^2 \ri ] = \de_k \bb{E} \li [  \li | A ( X - \mu ) \ri |^2 \ri ].
\ee   Substituting  (\ref{use88993b}) and (\ref{use88993c}) into
(\ref{use88993a}) yields  \be \la{keyfound}
 \bb{E} [ \mcal{T}_k ] \leq \tau  +
\de_k + \bb{E} \li [  \li | A ( X - \mu ) \ri |^2 \ri ] \leq \tau +
\de_k + \li | \li | A \ri | \ri |^2 \times \bb{E} [ || X - \mu ||^2
] \ee for all $k \in \bb{N}$.   Note that $\{ \mcal{T}_k, \; k \in
\bb{N} \}$ is a sequence of positive random variables. From Lemma
\ref{convergeas}, we have that $\mcal{T}_k \to \mcal{T}$ almost
surely as $k \to \iy$.   By Fatous' lemma, we have
\[
\liminf_{k \to \iy}  \bb{E} [\mcal{T}_k] \geq \bb{E} [ \liminf_{k
\to \iy} \mcal{T}_k] = \bb{E} [ \lim_{k \to \iy} \mcal{T}_k]  =
\bb{E} [\mcal{T}].
\]
Using (\ref{keyfound}), we have
\[
\liminf_{k \to \iy}  \bb{E} [\mcal{T}_k]  \leq \tau  + \bb{E} \li [
\li |  A ( X - \mu ) \ri |^2 \ri ] \leq \tau + \li | \li | A \ri |
\ri |^2 \times  \bb{E} \li [ \li | \li | X - \mu \ri | \ri |^2 \ri
].
\]
Hence, $\bb{E} [\mcal{T}]  \leq \tau  + \bb{E} \li [ \li | A (X -
\mu) \ri |^2 \ri ] \leq \tau  +   \li | \li | A \ri | \ri |^2 \times
\bb{E} \li [ \li | \li | X - \mu \ri | \ri |^2 \ri ]$.  Finally,
observing that $T \leq \mcal{T}$ as a consequence of the fact that
the continuity region $\mscr{R}$ of $T$ is a subset of the
continuity region $\mcal{R}$ of $\mcal{T}$, we have $\bb{E} [T] \leq
\bb{E} [\mcal{T}]$ and the desired results are proved.

\sect{Proof of Theorem \ref{Brown1} } \la{Brown1app}

From Theorem \ref{them8899Levy}, we know that $\bb{E} [ T ] < \iy$.
By Wald's equation for Brownian motion, we have that $\bb{E} [ W_{T}
]$ exists and  $\bb{E} [ W_{T} ] = \mu \bb{E} [ T ]$. Due to the
closedness of $\mscr{R}$ and  the continuity of the sample paths of
a Brownian motion, we have that $( T, W_{T} ) \in \mscr{R}$. By
Theorem \ref{fundamental}, we have $( \bb{E} [ T ], \bb{E} [ W_{T} ]
) \in \mscr{R}$.  Using Wald's equation, we have $( \bb{E} [ T ],
\mu \bb{E} [ T ] ) \in \mscr{R}$. Using this inclusion relation and
the definition of DET, we have $\bb{E} [ T ] \leq \mscr{B} (\mu)$.
This completes the proof of the theorem.

\sect{Proof of Corollary \ref{Brown3} } \la{Brown3app}

Note that $T = \{ t > 0: (t, W_t) \notin \mscr{R} \}$, where
$\mscr{R} = \{ (t, s): s \in \bb{R}^d, \; t > 0, \; t^2 \leq t g
(\f{s}{t} ) \}$. Since $g$ is a concave function on $\bb{R}^d$, it
follows from Lemma \ref{BoundaryCon} that $t g (\f{s}{t} )$ is a
concave function of $t > 0$ and $s \in \bb{R}^d$.  This implies that
the continuity region  $\mscr{R}$ is a convex set. Note that the DET
at $v = \mu$, is $\mscr{B} (\mu) = \sup \{ t
> 0: (t, \mu t) \in \mscr{R} \} = g (\mu) < \iy$.
It follows from Theorem \ref{Brown1} that $\bb{E} [ T ] \leq
\mscr{B} (\mu) = g (\mu)$. This completes the proof of the
corollary.

\sect{Proof of Theorem \ref{BoundConcentrationaLevya}}
\la{BoundConcentrationaLevyaapp}

Throughout the proof, all integrations are of Lebesgue type.  Note
that \be \la{Lintegral}
 \bb{E} [ T ]   =  \int_{\Om} T (\om) \; d \Pr \{ \om \}
  =  \int_{\Om} \li [ \int_{ \{ t: 0 < t < T (\om)    \} } d t \ri ]  d \Pr \{ \om
 \}.
\ee Applying Fubini's theorem to change the order of integration in
(\ref{Lintegral}) yields  \be
 \bb{E} [ T ]  =  \int_{0}^\iy \Pr \{ T > t \} dt =  \int_{0}^c  \Pr \li
\{ T > t \ri \} dt.   \la{VIP8899aLevy} \ee By the definition of $T$
and $\de (t)$, we have $\li \{ T > t \ri \} \subseteq \{ (t, X_t)
\in \mscr{R}  \} \subseteq  \li \{ || \ovl{X}_t - \mu || \geq \de
(t) \ri \}$ for $t > 0$.  It follows that \be \Pr \li \{ T
> t \ri \} \leq  \Pr \li \{ || \ovl{X}_t - \mu || \geq
\de (t) \ri \} \la{VIP8899bLevy} \ee for $t > 0$. Making use of
(\ref{VIP8899aLevy}) and (\ref{VIP8899bLevy}), we have $\bb{E} [ T ]
\leq \int_{0 }^c  \Pr \li \{  \li | \ovl{X}_t -  \mu \ri | \geq \de
(t) \ri  \} dt$.    This completes the proof of the theorem.

\end{document}